\documentclass{amsart}

\usepackage{amsmath,amssymb,amsthm}
\usepackage{array}
\usepackage{bbm,mathrsfs}
\usepackage{booktabs,tabularx}
\usepackage{float}

\newtheorem{theorem}{Theorem}[section]
\newtheorem*{theorem*}{Theorem}
\newtheorem*{maintheorem}{Main Theorem}
\newtheorem*{maincorollary}{Corollary}
\newtheorem{proposition}[theorem]{Proposition}
\newtheorem{corollary}[theorem]{Corollary}
\newtheorem{lemma}[theorem]{Lemma}
\theoremstyle{definition}
\newtheorem{definition}[theorem]{Definition}
\newtheorem{example}[theorem]{Example}
\newtheorem{assumption}[theorem]{Assumption}
\newtheorem{notation}{Notation}[section]
\newtheorem{remark}[theorem]{Remark}
\newtheorem*{remark*}{Remark}

\numberwithin{equation}{section}
\makeatletter
\renewcommand{\p@enumii}{\theenumi-}
\renewcommand{\p@enumiii}{\theenumi-\theenumii-}
\renewcommand{\p@enumiv}{\theenumi-\theenumii-\theenumiii-}
\makeatother

\newcommand{\hits}{\rightharpoonup}
\newcommand{\hitted}{\leftharpoonup}
\newcommand{\id}{\mathrm{id}}
\newcommand{\coo}{\mathrm{coinv}}
\newcommand{\ord}{\mathrm{ord}}
\renewcommand{\dim}{\mathrm{dim}}
\newcommand{\Ad}{\mathrm{Ad}}

\newcommand{\unit}{\pmb{1}}

\newcommand{\kk}{\Bbbk}
\newcommand{\G}{\mathbf{G}}
\newcommand{\Z}{\mathbb{Z}}

\newcommand{\DD}{\mathscr{D}}

\newcommand{\HH}{\mathcal{H}}
\renewcommand{\AA}{\mathcal{A}}
\newcommand{\LL}{\mathcal{L}}
\newcommand{\AKP}{A_{\mathrm{KP}}}

\newcommand{\ppp}{{\boldsymbol{p}}}
\newcommand{\ttt}{{\boldsymbol{t}}}
\newcommand{\lambdaa}{{\boldsymbol{\lambda}}}
\newcommand{\muu}{{\boldsymbol{\mu}}}

\newcommand{\ii}{\mathbbm{i}}
\newcommand{\fff}{{\boldsymbol{f}}}
\newcommand{\sss}{{\boldsymbol{s}}}

\newcommand{\sVect}{\mathcal{SV}}
\newcommand{\YD}[1]{{}^{#1}_{#1}\mathcal{YD}}

\newcommand{\AD}{\mathsf{AD}}
\newcommand{\SD}{\mathsf{SD}}
\newcommand{\CD}{\mathsf{CD}}

\newcommand{\Aut}{\mathrm{Aut}}

\newcommand{\inset}[2]{\left\{\begin{gathered}#1\end{gathered}\vphantom{\begin{gathered}#2\end{gathered}}\mathrel{}\middle|\mathrel{}\right. \left.\begin{gathered}#2\end{gathered}\vphantom{\begin{gathered}#1\end{gathered}}\right\}}
\newcommand{\bigangle}[2]{\left\langle\begin{gathered}#1\end{gathered}\vphantom{\begin{gathered}#2\end{gathered}}\mathrel{}\middle|\mathrel{}\right. \left.\begin{gathered}#2\end{gathered}\vphantom{\begin{gathered}#1\end{gathered}}\right\rangle}

\usepackage[hidelinks]{hyperref}

\begin{document}
\raggedbottom

\title[Hopf superalgebras of dimension up to $10$]{Completing the classification of Hopf superalgebras of dimension up to $10$}

\author[T.~Shibata]{Taiki Shibata}
\address[T.~Shibata]{Department of Applied Mathematics, Okayama University of Science, 1-1 Ridai-cho, Kita-ku, Okayama-shi, Okayama 700-0005, JAPAN}
\email{shibata@ous.ac.jp}
\author[K.~Shimizu]{Kenichi Shimizu}
\address[K.~Shimizu]{Department of Mathematical Sciences, Shibaura Institute of Technology, 307 Fukasaku, Minuma-ku, Saitama-shi, Saitama 337-8570, JAPAN}
\email{kshimizu@shibaura-it.ac.jp}
\author[R.~Wakao]{Ryota Wakao}
\address[R.~Wakao]{Office of the Mathematical Society of Japan, 1-34-8 Taito, Taito-ku, Tokyo 110-0016, JAPAN}
\email{msjadm01@air.ocn.ne.jp}

\date{\today}

\begin{abstract}
We introduce a family of semisimple Hopf superalgebras obtained from abelian extensions.
Using this family and bosonization, we complete the classification of Hopf superalgebras of dimension up to $10$ over an algebraically closed field of characteristic zero.
We also determine the duals of the semisimple classes with nonzero odd part and give tables of the complete classification.
\end{abstract}

\maketitle

\tableofcontents

\section{Introduction}\label{sec:intro}

The classification of finite-dimensional Hopf algebras of a prescribed dimension goes back to questions raised by Kaplansky and has been studied by many authors.
Among the results relevant to the present paper are the classifications in prime and prime-square dimensions~\cite{Zhu94,Mas96,AS98,Ng02}, in dimension $2p$~\cite{Mas95a,Ng05}, and in dimension $2p^2$~\cite{Nat99,HN09}, as well as classification results in dimensions $pq^2$ and $4p$~\cite{Nat99,Nat01,Nat04,CN11}.
The classifications and constructions in dimensions $8$, $12$, and $16$ are also used in this paper~\cite{KP66,Mas95b,Fuk97,Ste99,CDR00,Kas00,CDMM04,GV10}.
For counting methods and an overview of the classification problem, see~\cite{AN01,BG13}.

Throughout this paper, $\kk$ denotes an algebraically closed field of characteristic zero.
A Hopf superalgebra is a Hopf algebra object in the symmetric monoidal category $\sVect$ of supervector spaces over $\kk$.
Universal enveloping algebras and quantized enveloping algebras of Lie superalgebras are standard examples~\cite{Kac77,Yam94}.
The super setting also occurs in the classification of pointed Hopf algebras.
Nichols algebras and their liftings are central to that classification~\cite{AS98,AS10,AGI19a,AGI19b}, and Nichols algebras of diagonal type associated with contragredient Lie superalgebras were studied in~\cite{AAY11}.
More generally, Deligne's theorem shows that every finite symmetric tensor category over $\kk$ admits a faithful exact symmetric tensor functor to finite-dimensional supervector spaces~\cite{Del02}.

Some classes of Hopf superalgebras have general descriptions.
For example, the super version of the Cartier--Kostant theorem states that every finite-dimensional supercocommutative Hopf superalgebra is isomorphic to $\kk\Gamma\ltimes\bigwedge(V)$, where $\Gamma$ is a finite group and $V$ is a finite-dimensional purely odd $\kk\Gamma$-module~\cite{Kos77,AEG01}.
Taking duals gives the corresponding description of finite-dimensional supercommutative Hopf superalgebras.
Affine supergroups are represented by supercommutative Hopf superalgebras.
By~\cite{Mas05}, every supercommutative Hopf superalgebra is isomorphic as a superalgebra to the tensor product of its purely even quotient and an exterior algebra.
In the finitely generated case, these Hopf superalgebras are described by Harish--Chandra pairs~\cite{Mas12,MS18}.
Finite-dimensional Hopf superalgebras have also been studied using Heisenberg doubles and the graded pentagon equation~\cite{Dai24}.

The classification of Hopf superalgebras in fixed dimensions is less developed than its ordinary counterpart.
Aissaoui and Makhlouf~\cite{AM14} studied the cases of dimensions $2$, $3$, and $4$.
In~\cite{SW25}, we classified all Hopf superalgebras of dimension at most $5$ and constructed further classes in dimensions $6$ and $8$.
The non-semisimple pointed case in dimensions at most $10$ was treated in~\cite{SSW24}.
Throughout this paper, we use arXiv version~3 of both~\cite{SSW24} and~\cite{SW25}.
Taking duals gives further non-semisimple classes.
As in the ordinary case~\cite{Bea03}, the dual of a pointed Hopf superalgebra need not be pointed.
The remaining non-semisimple classes of dimension $8$, for which neither the Hopf superalgebra nor its dual is pointed, were obtained in~\cite{SW25} from the sixteen-dimensional Hopf algebras in~\cite{CDMM04,GV10}.

The purpose of this paper is to complete the classification in dimensions up to $10$, as announced at the end of~\cite{SSW24}.
The remaining cases are semisimple Hopf superalgebras with nonzero odd part.
They are necessarily non-pointed, since a finite-dimensional semisimple pointed Hopf superalgebra is purely even and isomorphic to a group algebra~\cite{SSW24}.
We introduce a family of semisimple Hopf superalgebras $\LL(\Gamma,\DD)$ of dimension $2|\Gamma|$.
Here $\Gamma$ is a finite group and $\DD=(\ttt,\ppp;\lambdaa,\muu)$ consists of an involutive automorphism $\ttt$ of $\Gamma$, a homomorphism $\ppp:\Gamma\to\Z_2=\Z/2\Z$, a normalized $1$-cochain $\lambdaa:\Gamma\to\kk^\times$, and a normalized $2$-cocycle $\muu:\Gamma\times\Gamma\to\kk^\times$, subject to the conditions in Definition~\ref{def:compatible-data}.
The structure maps are given in Section~\ref{sec:new}.
Theorems~\ref{thm:dim4-classification} and~\ref{thm:dim8-classification} show that all semisimple Hopf superalgebras with nonzero odd part in dimensions $4$ and $8$ belong to this family.

We use the Radford--Majid bosonization to obtain these Hopf superalgebras from ordinary Hopf algebras~\cite{Rad85,Maj94}.
The category $\sVect$ is a braided monoidal subcategory of the Yetter--Drinfeld category over $\kk C_2$, where $C_2$ is the cyclic group of order $2$.
The bosonization $\widehat{\HH}$ of a finite-dimensional Hopf superalgebra $\HH$ is therefore an ordinary Hopf algebra of dimension $2\dim(\HH)$.
In~\cite{SW25}, we gave a procedure for finding all $\HH$ with nonzero odd part such that $\widehat{\HH}\cong A$ for a fixed finite-dimensional Hopf algebra $A$.
It describes them by pairs $(g,\alpha)\in\G(A)\times\G(A^*)$, called super-data, modulo the action of the Hopf automorphism group of $A$.
Here $\G(A)$ denotes the group of group-like elements of $A$.
The associated Hopf superalgebra is obtained as a coinvariant subalgebra of $A$.
Thus ordinary classification results supply the possible bosonizations, while the super-data and their automorphism orbits determine the Hopf superalgebras.

For the semisimple cases, the ordinary classifications used here reduce the possible bosonizations to abelian extensions.
Natale's results give the possibilities in dimensions $2p^2$ and $4p$~\cite{Nat99,Nat04}; see also~\cite{CN11}.
The four-dimensional case uses the eight-dimensional Kac--Paljutkin Hopf algebra~\cite{KP66,Mas95b}, and the eight-dimensional case uses Kashina's classification in dimension $16$~\cite{Kas00}.
We determine the automorphism orbits of super-data for these extensions and compute the structure maps of the associated coinvariants.
For the cases arising in dimensions $4$ and $8$, these maps give the presentations $\LL(\Gamma,\DD)$ above.

Our results in dimensions $p^2$ and $2p$ hold for every odd prime $p$.
In dimension $2p$, we use the construction of Andruskiewitsch, Etingof, and Gelaki~\cite{AEG01}, recalled in Theorem~\ref{thm:AEG}.
It associates a Hopf superalgebra $A^u$ with a Hopf algebra $A$ and a group-like element $u\in\G(A)$ satisfying $u^2=1$.
For example, for an integer $m\geq3$, let
\[
 D_{2m}=\langle a,b\mid a^2=b^m=1,\ aba=b^{-1}\rangle
\]
be the dihedral group of order $2m$, where $a$ is a reflection and $b$ is a rotation.
The class occurring in dimension $2p$ is $(\kk D_{2p})^a$.

Our classification result is the following.

\begin{maintheorem}
The following statements hold.
\begin{enumerate}
\item For every odd prime $p$, every Hopf superalgebra of dimension $p^2$ is purely even.
\item For every odd prime $p$, there is exactly one semisimple non-pointed Hopf superalgebra of dimension $2p$ with nonzero odd part, up to isomorphism. It is represented by $(\kk D_{2p})^a$ and is self-dual.
\item In dimension $4$, there are exactly two semisimple non-pointed isomorphism classes with nonzero odd part, denoted by $\LL_4^{(1)}$ and $\LL_4^{(2)}$. Both are self-dual.
\item In dimension $8$, there are exactly twelve semisimple non-pointed isomorphism classes with nonzero odd part, denoted by $\LL_8^{(1)},\dots,\LL_8^{(12)}$.
\end{enumerate}
\end{maintheorem}

The dimension-$4$ classification in part~(3) was obtained in~\cite{SW25}; here we give compatible-data presentations of its two classes.
The classes $\LL_4^{(1)}$ and $\LL_4^{(2)}$ are defined in Theorem~\ref{thm:dim4-classification}.
The classes $\LL_8^{(1)},\dots,\LL_8^{(12)}$ are defined in Subsection~\ref{subsec:compatible-data} and classified in Theorem~\ref{thm:dim8-classification}.
We identify several of the eight-dimensional classes with familiar constructions, including $(\kk D_8)^a$ and the Hopf superalgebra obtained from the eight-dimensional Kac--Paljutkin algebra by the construction of Theorem~\ref{thm:AEG}, and we determine all duality relations among the twelve classes.

Combining the main theorem with~\cite{SSW24,SW25} and Appendix~\ref{subsec:app-overview} gives the following result.

\begin{maincorollary}
Up to isomorphism, there are exactly $92$ Hopf superalgebras of dimension at most $10$. Of these, $32$ are purely even, and the remaining $60$ are listed explicitly in Appendix~\ref{sec:appendix}.
\end{maincorollary}

Table~\ref{tab:sss-Hopf} gives the numbers of isomorphism classes of semisimple Hopf superalgebras in each dimension.
Representatives of the classes with nonzero odd part are listed in the last column.

\begin{table}[H]
\centering
\caption{Numbers of isomorphism classes of semisimple Hopf superalgebras of dimension $n\leq10$}
\label{tab:sss-Hopf}
\begin{tabular}{@{}c r r r l@{}}
\toprule
$n$ & Purely even & $\HH_{\bar1}\neq0$ & Total & Representatives \\
\midrule
1  & 1 & 0  & 1  & -- \\
2  & 1 & 0  & 1  & -- \\
3  & 1 & 0  & 1  & -- \\
4  & 2 & 2  & 4  & $\LL_4^{(1)},\LL_4^{(2)}$ \\
5  & 1 & 0  & 1  & -- \\
6  & 3 & 1  & 4  & $(\kk D_6)^a$ \\
7  & 1 & 0  & 1  & -- \\
8  & 8 & 12 & 20 & $\LL_8^{(1)},\dots,\LL_8^{(12)}$ \\
9  & 2 & 0  & 2  & -- \\
10 & 3 & 1  & 4  & $(\kk D_{10})^a$ \\
\bottomrule
\end{tabular}
\end{table}

\subsection*{Organization of the paper}
Section~\ref{sec:prelim} recalls the facts about Hopf superalgebras and bosonization used in the classification.
It explains how super-data and their equivalence classes describe the Hopf superalgebras with a fixed bosonization.
We also recall the construction $A^u$ associated with a group-like element $u$ satisfying $u^2=1$.
This construction is used both in the analysis of abelian extensions and in identifying several of the classes obtained later.

Section~\ref{sec:abel} develops the calculations needed for the semisimple cases.
We first describe abelian extensions and their Hopf automorphisms in terms of the finite groups and cochains defining the extension.
We then determine the structure of the coinvariant Hopf superalgebras associated with super-data.
These formulas lead to the family $\LL(\Gamma,\DD)$.
We prove that its members are semisimple and describe changes of the compatible datum that give isomorphic Hopf superalgebras.

Section~\ref{sec:class} applies these results to prove the main theorem.
We first treat dimensions $p^2$ and $2p$ for an arbitrary odd prime $p$.
We then determine the super-data up to automorphisms for the possible bosonizations in dimensions $8$ and $16$.
The resulting coinvariants give the compatible-data presentations of the four- and eight-dimensional classes.
In particular, the orbit calculations prove that these lists are complete and that their members are pairwise non-isomorphic.

Section~\ref{sec:dual} determines the duals of the eight-dimensional semisimple classes.
We use duality of the bosonizations together with explicit isomorphisms to distinguish the classes that share a bosonization.
Finally, Appendix~\ref{sec:appendix} combines these results with the non-semisimple classification.
It lists every class with nonzero odd part in dimensions at most $10$, together with its dual.
The counts also include the purely even classes, giving the corollary in the introduction.

\section{Preliminaries} \label{sec:prelim}
Throughout the paper, $\kk$ is an algebraically closed field of characteristic zero.
The unadorned symbol $\otimes$ means the tensor product over $\kk$.
In this section, we recall notation and results from~\cite{SSW24,SW25,AEG01}.
Some of these results hold over more general fields; see those papers for the precise assumptions.
We write $\Z_2=\{\bar0,\bar1\}$ for the additive group of order $2$.

\subsection{Hopf superalgebras} \label{subsec:Hopf-super}
A {\it supervector space} is a $\Z_2$-graded vector space $V=V_{\bar0}\oplus V_{\bar1}$, and $\sVect$ denotes the symmetric monoidal category of supervector spaces with grading-preserving linear maps.
Objects of $\sVect$ need not be finite-dimensional; its full subcategory of finite-dimensional objects is a symmetric tensor category.
A {\it Hopf superalgebra} is a Hopf algebra object in $\sVect$.
Here, a Hopf superalgebra $\HH$ is said to be {\it semisimple} (resp.~{\it pointed}) if
the category of left $\HH$-supermodules is semisimple (resp.~any simple right $\HH$-supercomodule is one-dimensional).
We focus on finite-dimensional semisimple non-pointed Hopf superalgebras.
Note that a finite-dimensional semisimple pointed Hopf superalgebra $\HH=\HH_{\bar0}\oplus\HH_{\bar1}$ is {\it purely even}
(i.e., the odd part $\HH_{\bar1}$ of $\HH$ is trivial)
and is isomorphic to a group algebra of a finite group.
See~\cite[Proposition~2.5]{SSW24}.

Recall that for a Hopf superalgebra $\HH$,
the bosonization $\widehat{\HH}:=\HH\otimes \kk C_2$ of $\HH$, where $C_2$ is the cyclic group of order $2$, becomes an ordinary Hopf algebra; see \cite[Section~2.2]{SSW24} for details.
We write $\G(\HH)$ for the group of even group-like elements of a Hopf superalgebra $\HH$.
The group $\G(\widehat{\HH})$ of group-like elements is isomorphic to $\G(\HH)\times C_2$.
Moreover, the category of left $\HH$-supermodules (respectively, supercomodules), with even morphisms, is monoidally equivalent to the category of left $\widehat{\HH}$-modules (respectively, comodules).
For finite-dimensional $\HH$, these equivalences restrict to tensor equivalences between the corresponding categories of finite-dimensional objects.
Thus $\HH$ is semisimple or pointed if and only if $\widehat{\HH}$ has the corresponding property.
The bosonization of a Hopf superalgebra having non-trivial odd part is neither commutative nor cocommutative.
Combining this fact with the classification of ordinary Hopf algebras of dimension $2p$~\cite{Mas95a,Ng05} gives the following result.
\begin{theorem}[\text{\cite[Theorem~5.2]{SW25}}]
All Hopf superalgebras of odd prime dimensions are purely even.
\end{theorem}

Let $\HH$ be a finite-dimensional Hopf superalgebra.
The dual space $\HH^*$ carries the Hopf superalgebra structure defined by the pairing convention in \cite[Section~3.4]{SW25}.
\begin{proposition}[\text{\cite[Proposition~2.6]{SSW24}}] \label{prp:dual-hat}
As a Hopf algebra, the bosonization $\widehat{\HH^*}$ of $\HH^*$ is isomorphic to the dual $\widehat{\HH}^*$ of $\widehat{\HH}$.
\end{proposition}

\subsection{Super-data and super-forms} \label{subsec:SD}
We recall the method used in~\cite{SSW24,SW25} to classify Hopf superalgebras.
For each finite-dimensional Hopf algebra $A$, the method determines all Hopf superalgebras $\HH$ such that $\widehat{\HH}\cong A$.
Write $C_2=\langle\omega\mid\omega^2=e\rangle$, and let $\YD{\kk C_2}$ denote the category of left--left Yetter--Drinfeld modules over $\kk C_2$.
We regard $\sVect$ as a full braided monoidal subcategory by giving each homogeneous vector $v$ the action $\omega\cdot v=(-1)^{|v|}v$ and the coaction $v\mapsto\omega^{|v|}\otimes v$.
We first consider Hopf algebras $\HH$ in $\YD{\kk C_2}$ such that $\widehat{\HH}\cong A$.

Let $A$ be a finite-dimensional Hopf algebra with comultiplication $\Delta$ and counit $\varepsilon$.
\begin{definition}
A pair $(g,\alpha) \in \G(A)\times \G(A^*)$ is called an {\it admissible datum} for $A$ if it satisfies
\[ g^2=1, \quad \alpha^2=\varepsilon \quad\text{and}\quad \alpha(g)=-1. \]
The set of all admissible data for $A$ is denoted by $\AD(A)$.
For $(g,\alpha),(g',\alpha')\in\AD(A)$, we write $(g,\alpha)\sim(g',\alpha')$ if there exists a Hopf algebra automorphism $\varphi$ of $A$ such that $\varphi(g)=g'$ and $\alpha=\alpha'\circ\varphi$. This is an equivalence relation on $\AD(A)$.
\end{definition}

For $(g,\alpha)\in\AD(A)$, the map
\[
\pi : A\longrightarrow \kk C_2; \quad
h\longmapsto \frac{\varepsilon(h)}{2}(e+\omega) + \frac{\alpha(h)}{2}(e-\omega)
\]
is a split epimorphism with section $\kk C_2\to A$ given by $\omega^i\mapsto g^i$.
The subalgebra of $\pi$-coinvariants of $A$ is
\[
A^{\coo(g,\alpha)} := \{ h \in A \mid (\id\otimes\pi)\Delta(h) = h \otimes e \}.
\]
The inverse procedure makes $A^{\coo(g,\alpha)}$ into a Hopf algebra in $\YD{\kk C_2}$ and yields the following result.
\begin{proposition}[\text{\cite[Proposition~4.11]{SW25}}] \label{prp:one-to-one}
Let $\mathscr{X}$ be the class of all Hopf algebras in $\YD{\kk C_2}$ whose bosonization is isomorphic to $A$ as a Hopf algebra.
The assignment $(g,\alpha)\mapsto A^{\coo(g,\alpha)}$ gives a bijection from $\AD(A) \mathord{/} \mathord{\sim}$ to $\mathscr{X} \mathord{/} \mathord{\cong}$.
\end{proposition}

\begin{definition}
We say that a Hopf superalgebra $\HH$ is a {\it super-form} of $A$ if $\widehat{\HH}$ is isomorphic to $A$.
If $\HH_{\bar1}\neq0$, then $\HH$ is called a {\it non-trivial super-form} of $A$.
\end{definition}

The following properties are proved in~\cite[Proposition~4.7]{SW25}.
\begin{proposition} \label{prp:super-forms}
Let $\HH$ be a super-form of $A$.
\begin{enumerate}
\item\label{prp:super-forms:(1)} Each of the groups $\G(A)$ and $\G(A^*)$ has a direct factor isomorphic to $C_2$.
\item\label{prp:super-forms:(2)} If the super-form $\HH$ is non-trivial, then $A$ is neither commutative nor cocommutative.
\item\label{prp:super-forms:(3)} If $A$ is semisimple/pointed, then so is $\HH$.
\end{enumerate}
\end{proposition}

By Proposition~\ref{prp:one-to-one}, the isomorphism classes of super-forms of $A$ are in bijection with the equivalence classes of $(g,\alpha)\in\AD(A)$ such that $A^{\coo(g,\alpha)}\in\YD{\kk C_2}$ belongs to $\sVect$.
By~\cite[Propositions~4.14 and~4.15]{SW25}, $\HH:=A^{\coo(g,\alpha)}$ is a super-form of $A$ if and only if
\[
\alpha\hits h \hitted\alpha = ghg
\quad\text{for all }h\in A,
\]
where $\alpha\hits h \hitted \alpha := (\alpha\otimes\id\otimes\alpha)(\id\otimes\Delta)\Delta(h)$.
When this condition holds, the $\Z_2$-grading of $\HH$ is
\[
\HH_{\bar\epsilon}=\{h\in\HH \mid ghg=(-1)^\epsilon h\}
\qquad(\epsilon\in\{0,1\}).
\]
In particular, the super-form $\HH$ is non-trivial if and only if $g\notin Z(A)$, where $Z(A)$ is the center of $A$.

\begin{definition} \label{def:SD}
An admissible datum $(g,\alpha)$ for $A$ is called a {\it super-datum} for $A$ if $g\notin Z(A)$ and $\alpha\hits h \hitted \alpha = ghg$ for all $h\in A$. The set of all super-data for $A$ is denoted by $\SD(A)$.
\end{definition}

By Proposition~\ref{prp:one-to-one}, the map $(g,\alpha) \mapsto A^{\coo(g,\alpha)}$ induces a one-to-one correspondence between $\SD(A)\mathord{/}\mathord{\sim}$ and the isomorphism classes of Hopf superalgebras $\HH$ such that $\HH_{\bar1}\neq0$ and $\widehat{\HH}\cong A$.

\subsection{Hopf superalgebras associated with involutions}\label{subsec:involution}
We recall the following construction from~\cite{AEG01}.

Let $A$ be a Hopf algebra with comultiplication $\Delta$ and counit $\varepsilon$, and take $u\in\G(A)$ with $u^2=1$.
As an algebra, put $A^u=A$, with grading
\[
A^u_{\bar\epsilon}:=\{h\in A\mid uhu=(-1)^\epsilon h\}
\qquad(\epsilon\in\{0,1\}).
\]
For homogeneous $h\in A^u$, write $\Delta(h)=\Delta_{\bar0}(h)+\Delta_{\bar1}(h)$ according to the parity of the second tensor factor.
\begin{theorem}[\text{\cite[Theorem~3.1.1]{AEG01}}] \label{thm:AEG}
The superalgebra $A^u$ becomes a Hopf superalgebra with
\[
\begin{gathered}
\Delta^u(h):=\Delta_{\bar0}(h)-(-1)^{|h|}(u\otimes1)\Delta_{\bar1}(h),
\quad \varepsilon^u(h):=\varepsilon(h),\\
S^u(h):=u^{|h|}S(h).
\end{gathered}
\]
This construction gives a bijective correspondence between pairs $(A,u)$ as above and pairs $(\HH,g)$ consisting of a Hopf superalgebra $\HH$ and an element $g\in\G(\HH)$ satisfying
$g^2=1$ and $gzg=(-1)^{|z|}z$ for all homogeneous $z\in\HH$.
\end{theorem}

For a finite group $G$ and an element $s\in G$ with $s^2=1$, we write $(\kk G)^s$ for the Hopf superalgebra obtained by taking $A=\kk G$ and $u=s$ in this construction.

\begin{example}
Let $G=D_6=\langle a,b\mid a^2=b^3=1,\ aba=b^{-1}\rangle$ and put
$c=(b+b^2)/2$ and $d=(b-b^2)/2$.
The Hopf superalgebra $(\kk D_6)^a$ has even basis $1,a,c,ac$ and odd basis $d,ad$.
Its multiplication is that of $\kk D_6$, and its comultiplication is determined by
\[
\Delta(a)=a\otimes a,\qquad
\Delta(c)=c\otimes c-ad\otimes d,\qquad
\Delta(d)=d\otimes c+ac\otimes d.
\]
For example, these formulas follow by applying Theorem~\ref{thm:AEG} to
$\Delta_{\kk D_6}(c)=c\otimes c+d\otimes d$ and
$\Delta_{\kk D_6}(d)=d\otimes c+c\otimes d$.
\end{example}

\section{Abelian extensions and their super-forms} \label{sec:abel}
In this section, we study the abelian extensions needed for our classification.
We recall their construction and describe their Hopf algebra automorphisms.
We then study their super-forms and introduce a family $\LL(\Gamma,\DD)$ of semisimple Hopf superalgebras.

\subsection{Abelian extensions}
For the general theory of Hopf algebra extensions and, in particular, abelian extensions, see, for example, \cite{Mas02}.

Let $G$ and $F$ be finite groups, and denote both identity elements by $1$.
Suppose that we are given a left action $\triangleright:G\times F\to F$ of $G$ on $F$ and a right action $\triangleleft:G\times F\to G$ of $F$ on $G$, written $x\triangleright u$ and $x\triangleleft u$, satisfying
\[
\begin{gathered}
x\triangleright (uv)=(x\triangleright u)((x\triangleleft u)\triangleright v),\\
(xy)\triangleleft u=(x\triangleleft(y\triangleright u))(y\triangleleft u)
\end{gathered}
\]
for all $x,y\in G$ and $u,v\in F$.
Then $(F,G;\triangleleft,\triangleright)$ is a {\it matched pair}.

Fix normalized $2$-cocycles $\sigma:F\times F\to(\kk^G)^\times$ and $\tau:G\times G\to(\kk^F)^\times$.
We identify them with the maps
\[
\begin{gathered}
G\times F\times F\longrightarrow\kk^\times,
\qquad (x,u,v)\longmapsto\sigma(x;u,v):=\sigma(u,v)(x),\\
G\times G\times F\longrightarrow\kk^\times,
\qquad (x,y,u)\longmapsto\tau(x,y;u):=\tau(x,y)(u).
\end{gathered}
\]
The map $\tau$ is called a dual $2$-cocycle.
Explicitly, for all $x,y,z\in G$ and $u,v,w\in F$, we require
\[
\begin{gathered}
\sigma(x\triangleleft u;v,w)\sigma(x;u,vw)
 =\sigma(x;uv,w)\sigma(x;u,v),\\
\sigma(x;1,u)=\sigma(x;u,1)=\sigma(1;u,v)=1,\\
\tau(x,y;z\triangleright u)\tau(xy,z;u)
 =\tau(x,yz;u)\tau(y,z;u),\\
\tau(x,1;u)=\tau(1,x;u)=\tau(x,y;1)=1.
\end{gathered}
\]
We also assume the compatibility condition
\[
\begin{aligned}
\sigma(xy;u,v)\tau(x,y;uv)
={}&\sigma\bigl(x;y\triangleright u,(y\triangleleft u)\triangleright v\bigr)
   \sigma(y;u,v)\\
&\times\tau(x,y;u)
   \tau\bigl(x\triangleleft(y\triangleright u),y\triangleleft u;v\bigr)
\end{aligned}
\]
for all $x,y\in G$ and $u,v\in F$.

Let $\kk^G:=(\kk G)^*=\bigoplus_{x\in G}\kk e_x$ be the dual Hopf algebra of the group algebra $\kk G$, where $e_x(y)=\delta_{x,y}$ for $x,y\in G$ and $\delta_{x,y}$ is the Kronecker delta.  Its unit is $\unit=\sum_{x\in G}e_x$.
The vector space
\[
\kk^G\#_{\sigma,\tau}\kk F:=\kk^G\otimes\kk F
 =\bigoplus_{x\in G,\,u\in F}\kk(e_x\otimes u)
\]
is a Hopf algebra with the following structure maps:
\begin{itemize}
\item (multiplication) $(e_x\otimes u)(e_y\otimes v):=\delta_{x\triangleleft u,y}\sigma(x;u,v)e_x\otimes uv$;
\item (unit) $\unit\otimes1=\sum_{x\in G}e_x\otimes1$;
\item (comultiplication) $\Delta(e_x\otimes u):=\sum_{yz=x}\tau(y,z;u)\bigl(e_y\otimes(z\triangleright u)\bigr)\otimes(e_z\otimes u)$;
\item (counit) $\varepsilon(e_x\otimes u):=\delta_{x,1}$.
\end{itemize}
The antipode is uniquely determined by these maps; see, for example, \cite{Mas02} for an explicit formula.
We abbreviate $e_x\otimes u$ by $e_xu$.
The Hopf algebra $A:=\kk^G\#_{\sigma,\tau}\kk F$ fits into an abelian extension
\[
1\longrightarrow\kk^G\longrightarrow A\longrightarrow\kk F\longrightarrow1
\]
and is called an {\it abelian extension} of $\kk F$ by $\kk^G$.

Set
\[
\begin{gathered}
\mathcal{G}(F,G):=
\inset{(\gamma,\fff)\in(\kk^G)^\times\times F}
{\forall x\in G,\ x\triangleright\fff=\fff,\\
 \forall x,y\in G,\ \gamma(xy)\tau(x,y;\fff)=\gamma(x)\gamma(y),\\
 \gamma(1)=1},\\
\mathcal{G}^*(F,G):=
\inset{(\sss,\beta)\in G\times(\kk^F)^\times}
{\forall u\in F,\ \sss\triangleleft u=\sss,\\
 \forall u,v\in F,\ \beta(uv)\sigma(\sss;u,v)=\beta(u)\beta(v),\\
 \beta(1)=1}.
\end{gathered}
\]
The maps
\begin{gather}
\mathcal{G}(F,G)\to\G(A);\quad (\gamma,\fff) \mapsto g_{\gamma,\fff} := \sum_{x\in G} \gamma(x) e_x \fff, \label{eq:G->G}\\
\mathcal{G}^*(F,G)\to\G(A^*);\quad (\sss,\beta) \mapsto (\alpha_{\sss,\beta} : e_xu \mapsto \delta_{x,\sss}\beta(u) ) \label{eq:G^*->G^*}
\end{gather}
are bijective.
Under this identification, $(g_{\gamma,\fff},\alpha_{\sss,\beta}) \in \AD(A)$ if and only if
\[
\begin{gathered}
\sss^2=1, \; \fff^2=1, \; \gamma(\sss)\beta(\fff)=-1,\\
\beta(u)\beta(\sss\triangleright u)\tau(\sss,\sss;u) =1
\quad\text{for all }u\in F,\\
\gamma(x)\gamma(x\triangleleft \fff) \sigma(x;\fff,\fff)=1
\quad\text{for all }x\in G.
\end{gathered}
\]

\subsection{Hopf algebra automorphisms of abelian extensions}
As in the previous subsection, let $A$ be an abelian extension $\kk^G\#_{\sigma,\tau}\kk F$ of $\kk F$ by $\kk^G$.
We describe the Hopf algebra automorphisms of $A$ under the following assumptions.
\begin{assumption}\label{ass:actions}
We assume the following two conditions.
\begin{enumerate}
\item \label{assump:auto1} The group $G$ is abelian and $\G(A)\cong G$ as groups.
\item \label{assump:auto2} The group $F$ is $C_2=\langle t\rangle$, and the left action $\triangleright:G\times F\to F$ is trivial.
\end{enumerate}
\end{assumption}
Since $\kk$ is algebraically closed of characteristic zero, the assumption that $G$ is abelian implies
$\kk^G\cong\kk G$ as Hopf algebras.
Thus Assumption~\ref{ass:actions}\eqref{assump:auto1} gives
\[
\kk^G=\kk[\G(A)].
\]

Let $\varphi$ be a Hopf algebra automorphism of $A$.
The equality above implies that $\varphi$ preserves $\kk^G$.
Hence there is a group automorphism $\overline{\varphi}$ of $G$ such that
\[
\varphi(e_x)=e_{\overline{\varphi}(x)}
\quad\text{for all }x\in G.
\]
Let $(\kk^G)^+$ denote the kernel of the counit of $\kk^G$.  The induced automorphism of the quotient $A/(\kk^G)^+A\cong\kk C_2$ is the identity, since $C_2$ has no non-trivial automorphism.
Writing $q:A\to\kk C_2$ for the quotient map, we therefore have $q\circ\varphi=q$.
Thus $\varphi$ preserves the right coaction $\rho=(\id\otimes q)\Delta$.
Since $\rho(e_xt^\ell)=e_xt^\ell\otimes t^\ell$ for $\ell\in\{0,1\}$, $\varphi$ preserves the subspace $\kk^Gt$.
Combining this with $\varphi(e_x)=e_{\overline{\varphi}(x)}$, we obtain scalars $c_x\in\kk^\times$ such that
\[
\varphi(e_xt)=c_xe_{\overline{\varphi}(x)}t
\quad\text{for all }x\in G.
\]
The counit gives $c_1=1$.

The multiplication in $A$ gives
\[
(e_xt)e_y=\delta_{x\triangleleft t,y}e_xt.
\]
Applying $\varphi$ to this equality shows that
\[
\overline{\varphi}(x\triangleleft t)=\overline{\varphi}(x)\triangleleft t
\quad\text{for all }x\in G.
\]
We identify $t\in C_2$ with $\unit\#_{\sigma,\tau}t\in A$.  In $A$,
\[
t^2=(\unit\#_{\sigma,\tau}t)^2
   =\sum_{x\in G}\sigma(x;t,t)e_x.
\]
Applying $\varphi$ to this equality gives
\[
c_xc_{x\triangleleft t}
 =\sigma(x;t,t)\sigma(\overline{\varphi}(x);t,t)^{-1}
\quad\text{for all }x\in G.
\]
Finally, comparing $\Delta(\varphi(e_xt))$ with
$(\varphi\otimes\varphi)\Delta(e_xt)$ gives
\[
c_xc_y
 =c_{xy}\tau(x,y;t)^{-1}
  \tau(\overline{\varphi}(x),\overline{\varphi}(y);t)
\quad\text{for all }x,y\in G.
\]

\begin{proposition} \label{prp:automorphism}
For a Hopf algebra automorphism $\varphi$ of $A$, there exist a group automorphism $\overline{\varphi}$ of $G$ and scalars $c_x\in\kk^\times$ ($x\in G$) such that
\begin{equation}\label{eq:automorphism-form}
\varphi(e_xt^\ell)=c_x^\ell e_{\overline{\varphi}(x)}t^\ell
\quad\text{for all }x\in G\text{ and }\ell\in\{0,1\}.
\end{equation}
Moreover,
\begin{equation} \label{eq:action-automorphism}
\overline{\varphi}(x\triangleleft t)=\overline{\varphi}(x)\triangleleft t
\quad\text{for all }x\in G,
\end{equation}
and
\begin{equation} \label{eq:automorphism}
\begin{gathered}
c_1=1,\\
c_xc_{x\triangleleft t}
 =\sigma(x;t,t)\sigma(\overline{\varphi}(x);t,t)^{-1}
 \quad\text{for all }x\in G,\\
c_xc_y
 =c_{xy}\tau(x,y;t)^{-1}
  \tau(\overline{\varphi}(x),\overline{\varphi}(y);t)
 \quad\text{for all }x,y\in G.
\end{gathered}
\end{equation}
Conversely, every pair $(\overline{\varphi},(c_x)_{x\in G})$ satisfying
\eqref{eq:action-automorphism} and \eqref{eq:automorphism} defines a Hopf algebra automorphism by \eqref{eq:automorphism-form}.
\end{proposition}
\begin{proof}
The necessity of the conditions was proved above.
Conversely, the conditions ensure that the map preserves multiplication and comultiplication on the basis elements, as well as the unit and counit.
The resulting bijective bialgebra map automatically commutes with the antipode.
\end{proof}

\begin{lemma}\label{lem:automorphism-count}
Let $m$ be the number of group automorphisms $\psi$ of $G$ satisfying
\[
\psi(x\triangleleft t)=\psi(x)\triangleleft t,
\qquad
\frac{\tau(\psi(x),\psi(y);t)}{\tau(\psi(y),\psi(x);t)}
=\frac{\tau(x,y;t)}{\tau(y,x;t)}
\quad(x,y\in G),
\]
and let $n$ be the number of characters $\chi:G\to\kk^\times$ satisfying
\[
\chi(x)\chi(x\triangleleft t)=1\quad(x\in G).
\]
Then $|\Aut(A)|\leq mn$.
For each fixed $\overline\varphi\in\Aut(G)$, the number of solutions of~\eqref{eq:automorphism} is either $0$ or $n$.
\end{lemma}
\begin{proof}
By Proposition~\ref{prp:automorphism}, the group automorphism induced by a Hopf automorphism commutes with the action of $t$.
Interchanging $x$ and $y$ in the last equation of~\eqref{eq:automorphism} gives the stated condition on $\tau$, since $G$ is abelian.
For a fixed $\overline\varphi$, the pointwise ratio of two solutions of~\eqref{eq:automorphism} is a character satisfying the stated condition.
Conversely, multiplying a solution pointwise by such a character gives another solution.
Thus each of the $m$ possible group automorphisms has at most $n$ lifts to $A$.
\end{proof}

Fix $\alpha,\alpha'\in\G(A^*)$ and take the corresponding pairs
$(\sss,\beta),(\sss',\beta')\in \mathcal{G}^*(F,G)$ under the bijection~\eqref{eq:G^*->G^*}, respectively.
Let $\varphi$ be a Hopf algebra automorphism of $A$.
Then Proposition~\ref{prp:automorphism} gives a group automorphism $\overline{\varphi}$ of $G$ and scalars $c_x\in\kk^\times$ satisfying \eqref{eq:automorphism-form}.

\begin{proposition} \label{prp:automorphism2}
$\alpha=\alpha'\circ \varphi$ if and only if
$\overline{\varphi}(\sss)=\sss'$ and $\beta(t)=c_\sss\beta'(t)$.
\end{proposition}

\subsection{Super-forms of abelian extensions} \label{sec:ell=0}
In this subsection, we describe the super-forms of an abelian extension $A:=\kk^G\#_{\sigma,\tau}\kk F$ of $\kk F$ by $\kk^G$
satisfying Assumption~\ref{ass:actions}\eqref{assump:auto2},
that is, the left action $\triangleright$ is trivial and $F=C_2=\langle t\rangle$.
Assumption~\ref{ass:actions}\eqref{assump:auto1} is needed later for the calculation of automorphism groups, but it is not assumed in this subsection.

We take $(\gamma,t^\ell)\in\mathcal{G}(F,G)$ with $\ell\in\{0,1\}$ and $(\sss,\beta)\in\mathcal{G}^*(F,G)$.
Then the pair $(g_{\gamma,t^\ell},\alpha_{\sss,\beta})$ is an admissible datum for $A$ if and only if
\[
\begin{aligned}
\sss^2&=1,\\
\gamma(x)\gamma(x\triangleleft t^\ell)
 &=\sigma(x;t^\ell,t^\ell)^{-1}
 &&\text{for all }x\in G,\\
\gamma(\sss)&=-\beta(t^\ell)^{-1},\\
\beta(t)^2&=\tau(\sss,\sss;t)^{-1}.
\end{aligned}
\]
Fix $(g,\alpha)=(g_{\gamma,t^\ell},\alpha_{\sss,\beta})\in\AD(A)$.
A direct computation shows the following.
\begin{lemma} \label{prp:SD(abel)}
The pair $(g,\alpha)$ belongs to $\SD(A)$ if and only if the following conditions hold:
\begin{enumerate}
\item[(i)] $x\triangleleft t^\ell=\Ad_\sss(x)$ for all $x\in G$;
\item[(ii)] $\gamma(x)^{-1}\gamma(x\triangleleft t)=\tau(\sss,\Ad_\sss^\ell(x);t)\tau(x,\sss;t)^{-1}$ for all $x\in G$;
\item[(iii)] $\Ad_\sss^\ell(w)\neq w$ or $\gamma(w)\neq\gamma(w\triangleleft t)$ for some $w\in G$.
\end{enumerate}
Here $\Ad_\sss:G\to G;\ x\mapsto\sss x\sss$.
\end{lemma}
\begin{proof}
Apply the identity $\alpha\hits h\hitted\alpha=ghg$ first to $h=e_x$ and then to $h=e_xt$.
Comparison of the indices for $h=e_x$ gives (i), and comparison of the coefficients for $h=e_xt$ gives (ii).
Condition (iii) says that $g$ fails to commute with at least one of these basis elements, or equivalently that $g\notin Z(A)$.
\end{proof}

For the rest of this subsection, assume that $(g,\alpha)\in\SD(A)$.
We describe the super-form $\HH=A^{\coo(g,\alpha)}$ for $\ell=0,1$.

\subsubsection{The case \texorpdfstring{$\ell=1$}{ell=1}}
Suppose that $\ell=1$ and put $u=g$.
Since $u=\sum_x\gamma(x)e_xt$ is a group-like element with $u^2=1$,
replacing $t$ by $u$ makes both cocycles trivial.
By Lemma~\ref{prp:SD(abel)}, $\sss$ is a non-central involution in $G$
and $x\triangleleft t=\Ad_\sss(x)$ for all $x\in G$.
In particular, this case does not occur when $G$ is abelian.
Thus
\[
A=\kk^G\rtimes\kk\langle u\rangle,
\qquad u^2=1,\qquad
u f u=f\circ\Ad_\sss\quad(f\in\kk^G).
\]
Moreover, $\alpha(e_x)=\delta_{x,\sss}$ and $\alpha(u)=-1$, so
$\alpha(e_xu^j)=\delta_{x,\sss}(-1)^j$ for $x\in G$ and $j\in\{0,1\}$.
We now describe the corresponding super-form.

\begin{proposition}\label{prp:ell-one-superform}
The super-form $\HH=A^{\coo(g,\alpha)}$ is isomorphic to the following
Hopf superalgebra on the vector space $\kk^G$:
\[
\HH_{\bar\epsilon}=\{f\in\kk^G\mid f(\sss x\sss)=(-1)^\epsilon f(x)
\text{ for all }x\in G\}\quad(\epsilon\in\{0,1\}).
\]
The coalgebra structure is the usual one on $\kk^G$, dual to the group algebra $\kk G$.
The product is given by
\[
(f\star h)(x)=f(x\sss^{|h|})h(x)
\]
for homogeneous $h$.
The unit is the constant function $1$, and the antipode is given by
$S(f)(x)=f(\sss^{|f|}x^{-1})$ for homogeneous $f$.
\end{proposition}
\begin{proof}
Put $R_\sss f(x)=f(x\sss)$ and $\vartheta f=f\circ\Ad_\sss$.
The map $(\id\otimes\alpha)\Delta$ sends
$f u^j$ to $(-1)^j(R_\sss f)u^j$.
Consequently, the linear map
\[
J:\kk^G\longrightarrow\HH;\quad
f\longmapsto\frac12(f+R_\sss f)+\frac12(f-R_\sss f)u
\]
is a bijection.
Moreover, $uJ(f)u=J(\vartheta f)$, so it preserves the stated grading.
For homogeneous $h$, multiplication in $A$ gives
\[
J(f)J(h)=J\bigl((R_\sss^{|h|}f)h\bigr).
\]
Write $\Delta_0$ for the ordinary comultiplication of $\kk^G$.
The identities
\[
\Delta_0 R_\sss=(\id\otimes R_\sss)\Delta_0,
\qquad
(R_\sss\otimes\id)\Delta_0
=(\id\otimes R_\sss\vartheta)\Delta_0
\]
show that, after decomposing the tensor factors into homogeneous terms,
\[
\Delta_A(J(f))
=\sum J(f_{(1)})u^{|f_{(2)}|}\otimes J(f_{(2)}),
\qquad
\Delta_0(f)=\sum f_{(1)}\otimes f_{(2)}.
\]
Thus the inverse bosonization gives the stated comultiplication and counit.
The antipode formula follows from
$S_{\HH}(J(f))=u^{|f|}S_A(J(f))$ for homogeneous $f$.

\end{proof}

\begin{corollary}\label{cor:ell-one-identification}
The super-form in Proposition~\ref{prp:ell-one-superform} satisfies
\[
\HH\cong\bigl((\kk G)^{\sss}\bigr)^*,
\]
where $(-)^{\sss}$ denotes the construction in Theorem~\ref{thm:AEG}.
\end{corollary}
\begin{proof}
Use the evaluation pairing $\langle x,f\rangle=f(x)$
and the tensor-pairing convention fixed in Subsection~\ref{subsec:Hopf-super}.
The dual has underlying algebra $\kk G$, grading given by
$z\mapsto\sss z\sss$, and comultiplication
\[
\Delta(x)=x\otimes x_++x\sss\otimes x_-,
\qquad
x_\pm=\frac12(x\pm\sss x\sss)
\quad(x\in G).
\]
For homogeneous $z$, this is
$\Delta_{\bar0}(z)-(-1)^{|z|}(\sss\otimes1)\Delta_{\bar1}(z)$,
where the subscripts refer to the parity of the second tensor factor
of the ordinary group-algebra comultiplication.
This is precisely the structure of $(\kk G)^{\sss}$.
\end{proof}

\subsubsection{The case \texorpdfstring{$\ell=0$}{ell=0}}
Suppose that $\ell=0$. Then
\begin{equation}\label{eq:t^ell}
g=g_{\gamma,1}=\sum_{x\in G}\gamma(x)e_x.
\end{equation}
Lemma~\ref{prp:SD(abel)} says that $\sss$ is central and that
\begin{equation}\label{eq:SD(abel)}
\begin{aligned}
\gamma(x)\gamma(x\triangleleft t)
 &=\tau(\sss,x;t)\tau(x,\sss;t)^{-1}
 &&\text{for all }x\in G,\\
\gamma(w)&\neq\gamma(w\triangleleft t)
 &&\text{for some }w\in G.
\end{aligned}
\end{equation}
Since $\gamma:G\to\{\pm1\}$ is a group homomorphism
and $\gamma(\sss)=-1$, its kernel has index $2$.
Hence, on putting $\Gamma:=\operatorname{Ker}(\gamma)$, we have
$G=\Gamma\times\langle\sss\rangle$.
For $x\in\Gamma$, define $\ttt(x)\in\Gamma$ and
$\ppp(x)\in\{0,1\}$ by $x\triangleleft t=\ttt(x)\sss^{\ppp(x)}$.
Set
\begin{equation}\label{eq:induced-mu}
\begin{aligned}
c_x&:=\bigl(\beta(t)\tau(x,\sss;t)\bigr)^{-1},&
\lambdaa(x)&:=c_{\ttt(x)}^{\ppp(x)}\sigma(x;t,t),\\
\muu(x,y)&:=(-1)^{\ppp(x)\ppp(y)}\tau(x,y;t),&
B_{x,\bar\epsilon}&:=(e_x+c_x^\epsilon e_{x\sss})t^\epsilon
\end{aligned}
\end{equation}
for $x,y\in\Gamma$ and $\epsilon\in\{0,1\}$.

\begin{proposition}\label{prp:ell0-normal-form}
The map $\ttt$ is an involutive group automorphism of $\Gamma$,
and $\ppp:\Gamma\to\Z_2$ is a non-trivial group homomorphism
satisfying $\ppp\circ\ttt=\ppp$.
The elements $B_{x,\bar\epsilon}$ form a homogeneous basis of
$\HH=A^{\coo(g,\alpha)}$, with
$|B_{x,\bar\epsilon}|=\epsilon\ppp(x)$.
The structure maps are
\[
\begin{aligned}
B_{x,\bar\epsilon}B_{y,\bar\eta}
 &=\delta_{\ttt^\epsilon(x),y}\lambdaa(x)^{\epsilon\eta}
   B_{x,\overline{\epsilon+\eta}},\\
\Delta_{\HH}(B_{x,\bar\epsilon})
 &=\sum_{yz=x}\muu(y,z)^\epsilon
   B_{y,\bar\epsilon}\otimes B_{z,\bar\epsilon},\\
1_{\HH}&=\sum_{x\in\Gamma}B_{x,\bar0},
\qquad
\varepsilon_{\HH}(B_{x,\bar\epsilon})=\delta_{x,1}.
\end{aligned}
\]
Here the products in the sum are taken in $\Gamma$.
\end{proposition}

\begin{proof}
The action $(-)\triangleleft t$ is an involutive group automorphism
of $G$ fixing $\sss$.
Applying it to the product $xy$ and to $x\triangleleft t$ shows
that $\ttt$ and $\ppp$ are group homomorphisms,
$\ttt^2=\id$, and $\ppp\circ\ttt=\ppp$.
Since $\gamma(x)=1$ for $x\in\Gamma$,
Lemma~\ref{prp:SD(abel)} gives
\begin{equation}\label{eq:ell0-tau-sign}
\tau(\sss,x;t)=(-1)^{\ppp(x)}\tau(x,\sss;t).
\end{equation}
If $\ppp$ were trivial, $\gamma$ would be invariant under
$(-)\triangleleft t$, making $g$ central in $A$.
Thus $\ppp$ is non-trivial.

Write $\tau(a,b)=\tau(a,b;t)$ and $\sigma(a)=\sigma(a;t,t)$
within this proof, and put $E_x=B_{x,\bar0}$ and
$F_x=B_{x,\bar1}$.
The coinvariants are the fixed points of
$\rho=(\id\otimes\alpha)\Delta_A$.
The cocycle identity and $\beta(t)^2\tau(\sss,\sss)=1$ give
\[
\rho(e_x)=e_{x\sss},\qquad
\rho(e_xt)=c_xe_{x\sss}t,\qquad
\rho(e_{x\sss}t)=c_x^{-1}e_xt.
\]
Therefore the $E_x,F_x$ form a basis of $\HH$.
Conjugation by $g$ fixes $E_x$ and sends $F_x$ to
$(-1)^{\ppp(x)}F_x$, proving the parity assertion.

The products involving an $E_x$ follow from the orthogonal
idempotents $e_z$ and the defining formula for $\ttt$.
For the remaining product, compatibility of $\sigma$ and $\tau$
gives
\[
\sigma(x\sss)=\sigma(x)\beta(t)^2\tau(x,\sss)\tau(x\triangleleft t,\sss)
=\sigma(x)c_x^{-1}c_{\ttt(x)}^{\,2\ppp(x)-1}.
\]
Here the second equality uses
$\tau(z\sss,\sss)=\tau(\sss,\sss)/\tau(z,\sss)$.
The product $F_xF_y$ vanishes unless $y=\ttt(x)$.
In that case its coefficient at $e_x$ is
$c_{\ttt(x)}^{\ppp(x)}\sigma(x)=\lambdaa(x)$,
and this identity gives the same coefficient at
$e_{x\sss}$.
This proves the multiplication formula.

The coproduct of $E_z$ is
$\Delta_A(E_z)=\sum_{xy=z}E_x\otimes E_y$.
To compute that of $F_z$, the cocycle identity and
\eqref{eq:ell0-tau-sign} give, for $x,y\in\Gamma$,
\[
\begin{aligned}
c_{xy}\tau(x,y\sss)&=\tau(x,y)c_y,\\
c_{xy}\tau(x\sss,y)&=(-1)^{\ppp(y)}\tau(x,y)c_x,\\
\tau(x\sss,y\sss)&=(-1)^{\ppp(y)}\tau(x,y)c_xc_y.
\end{aligned}
\]
Moreover,
\[
F_xg^{\ppp(y)}
 =(-1)^{\ppp(x)\ppp(y)}
   (e_x+(-1)^{\ppp(y)}c_xe_{x\sss})t.
\]
Comparing the coefficients of
$e_xt\otimes e_yt$, $e_xt\otimes e_{y\sss}t$,
$e_{x\sss}t\otimes e_yt$, and
$e_{x\sss}t\otimes e_{y\sss}t$ now yields
\[
\Delta_A(F_z)
 =\sum_{xy=z}\muu(x,y)
   F_xg^{\ppp(y)}\otimes F_y.
\]
Since $|F_y|=\ppp(y)$, the bosonization formula gives
$\Delta_{\HH}(F_z)=\sum_{xy=z}\muu(x,y)F_x\otimes F_y$.
This proves the coproduct formula.
The unit and counit formulas follow from their formulas in $A$.
\end{proof}

Since $\ppp:\Gamma\to\Z_2$ is non-trivial, $|\Gamma|$ is even and $\dim(\HH)=2|\Gamma|$ is divisible by $4$.
Thus a super-form of dimension $2p$, with $p$ an odd prime, can occur only in the case $\ell=1$.

The restriction of the quotient map $G\to G/\langle\sss\rangle$ to
$\Gamma=\operatorname{Ker}(\gamma)$ is a group isomorphism.
In the next subsection, we use these formulas to define the family $\LL(\Gamma,\DD)$.

\subsection{The family \texorpdfstring{$\LL(\Gamma,\DD)$}{L(Gamma,D)}} \label{sec:new}
Let $\Gamma$ be a finite group.
We use standard group-cohomology notation.
For $n\in\mathbb{N}$, we set $C^n(\Gamma):=\{ f:\Gamma^n\to\kk^\times \mid f(x_1,\dots,x_n)=1 \text{ if } \exists i \text{ s.t. } x_i=1 \}$
and $Z^n(\Gamma):=\operatorname{Ker}(\partial_n)$,
where $\partial_n:C^n(\Gamma)\to C^{n+1}(\Gamma)$ is the $n$th coboundary map.
For example, $\partial_1 f(x,y) = f(x)f(y)f(xy)^{-1}$ for $f\in C^1(\Gamma)$.
Let $H^n(\Gamma):=Z^n(\Gamma)/\operatorname{Im}(\partial_{n-1})$ be the $n$th cohomology group.

\begin{definition}\label{def:compatible-data}
Suppose that we are given the following maps.
\begin{enumerate}
\item An involutive group automorphism $\ttt$ of $\Gamma$.
\item A map $\ppp:\Gamma \to \{0,1\}$ such that $\ppp$ is a homomorphism under the identification $\{0,1\}=\mathbb{Z}_2$ and $\ppp\circ\ttt=\ppp$.
\item $\lambdaa \in C^1(\Gamma)$ satisfying $\lambdaa\circ\ttt=\lambdaa$.
\item $\muu\in Z^2(\Gamma)$.
\end{enumerate}
We say that the quadruple $(\ttt,\ppp;\lambdaa,\muu)$ is a {\it compatible datum} for $\Gamma$ if the following equation holds
for all $x,y\in\Gamma$.
\begin{equation} \label{eq:compatible}
\lambdaa(xy) = (-1)^{\ppp(x)\ppp(y)} \lambdaa(x)\lambdaa(y)\muu(x,y)\muu(\ttt(x),\ttt(y)).
\end{equation}
\end{definition}

Let $\DD=(\ttt,\ppp;\lambdaa,\muu)$ be a compatible datum for $\Gamma$
and let the $2|\Gamma|$-dimensional vector space
\[
\LL(\Gamma,\DD) := \bigoplus_{x\in\Gamma, \epsilon\in\{0,1\}} \kk v(x,\bar{\epsilon})
\]
be a superspace by setting $|v(x,\bar{\epsilon})| := \delta_{\epsilon,1} \ppp(x)$.
Define structure maps by
\begin{itemize}
\item $v(x,\bar{\epsilon})v(y,\bar{\eta}) := \delta_{\ttt^\epsilon(x),y} \lambdaa(x)^{\epsilon \eta} v(x,\overline{\epsilon+\eta})$;
\item $1=\sum_{x\in\Gamma} v(x,\bar{0})$;
\item $\Delta(v(x,\bar{\epsilon})) := \sum_{yz=x} \muu(y,z)^\epsilon v(y,\bar{\epsilon}) \otimes v(z,\bar{\epsilon})$;
\item $\varepsilon(v(x,\bar{\epsilon})) := \delta_{x,1}$;
\item $S(v(x,\bar{\epsilon})) := (\lambdaa(x^{-1})\muu(x^{-1},x))^{-\epsilon} v(\ttt^\epsilon(x^{-1}),\bar{\epsilon})$.
\end{itemize}
Here, $x,y\in\Gamma$ and $\epsilon,\eta\in\{0,1\}$.

\begin{proposition}\label{prp:compatible-Hopf}
The preceding formulas make $\LL(\Gamma,\DD)$ into a Hopf superalgebra.
\end{proposition}
\begin{proof}
Associativity follows from $\ttt^2=\id$ and $\lambdaa\circ\ttt=\lambdaa$, while coassociativity is the $2$-cocycle identity for $\muu$.
We check the multiplicativity of $\Delta$ on basis elements.
Both sides vanish unless the Kronecker-delta condition in the product is satisfied, and the cases $(\epsilon,\eta)\neq(1,1)$ follow directly from the definitions.
For the remaining case, the multiplication in the tensor product of superspaces gives
\[
\begin{aligned}
&\Delta(v(x,\bar1))\Delta(v(\ttt(x),\bar1))\\
&\quad=
\sum_{uv=x}
(-1)^{\ppp(u)\ppp(v)}
\lambdaa(u)\lambdaa(v)
\muu(u,v)\muu(\ttt(u),\ttt(v))
\,v(u,\bar0)\otimes v(v,\bar0)\\
&\quad=
\lambdaa(x)\sum_{uv=x}v(u,\bar0)\otimes v(v,\bar0)\\
&\quad=
\Delta\bigl(v(x,\bar1)v(\ttt(x),\bar1)\bigr),
\end{aligned}
\]
where the second equality is exactly \eqref{eq:compatible}.
The unit and counit axioms are immediate.
For the antipode, the cocycle identity gives
$\muu(x,x^{-1})=\muu(x^{-1},x)$.
Using this identity and $\lambdaa\circ\ttt=\lambdaa$, substitution gives
\[
m(S\otimes\id)\Delta(v(x,\bar\epsilon))
=m(\id\otimes S)\Delta(v(x,\bar\epsilon))
=\delta_{x,1}\sum_{z\in\Gamma}v(z,\bar0).
\]
Thus $S$ is an antipode.
\end{proof}

We denote by $\CD(\Gamma)$ the set of compatible data $\DD=(\ttt,\ppp;\lambdaa,\muu)$ for which $\ppp$ is non-trivial; equivalently, $\LL(\Gamma,\DD)_{\bar1}\neq0$.

\begin{proposition} \label{prp:induced-compatible-datum}
For the case $\ell=0$ of Section~\ref{sec:ell=0}, put
$\Gamma=\operatorname{Ker}(\gamma)$ and take the maps
$\ttt,\ppp,\lambdaa,\muu$ defined there.
Then $\DD=(\ttt,\ppp;\lambdaa,\muu)$ belongs to $\CD(\Gamma)$, and
\[
B_{x,\bar\epsilon}\longmapsto v(x,\bar\epsilon)
\]
extends to an isomorphism of Hopf superalgebras
$A^{\coo(g,\alpha)}\cong\LL(\Gamma,\DD)$.
\end{proposition}
\begin{proof}
Proposition~\ref{prp:ell0-normal-form} proves that $\ttt$ is an involutive
group automorphism and that $\ppp$ is a non-trivial homomorphism with
$\ppp\circ\ttt=\ppp$.
Associativity of $B_{x,\bar1}B_{\ttt(x),\bar1}B_{x,\bar1}$ gives
$\lambdaa\circ\ttt=\lambdaa$.
The formula for $\muu$ is the product of the restriction of the normalized
$2$-cocycle $\tau(-,-;t)$ and the $2$-cocycle
$(x,y)\mapsto(-1)^{\ppp(x)\ppp(y)}$; hence $\muu\in Z^2(\Gamma)$.
Normalization also gives $\lambdaa(1)=1$.
Finally, comparing the coefficient of
$B_{x,\bar0}\otimes B_{y,\bar0}$ in
$\Delta(B_{xy,\bar1}B_{\ttt(xy),\bar1})
=\Delta(B_{xy,\bar1})\Delta(B_{\ttt(xy),\bar1})$
gives~\eqref{eq:compatible}.
The structure maps in Proposition~\ref{prp:ell0-normal-form} now give
the asserted isomorphism.
\end{proof}
The super-datum determines the group $\Gamma=\operatorname{Ker}(\gamma)$.
In the examples below, we identify it with a standard group and, when useful,
rescale the basis $B_{x,\bar1}$.
These rescalings are described by equivalent compatible data, as defined below.

For every super-datum used below in dimensions $4$ and $8$, one has $\fff=1$ (equivalently, $\ell=0$ in Section~\ref{sec:ell=0}).
Thus Proposition~\ref{prp:induced-compatible-datum} applies to every non-trivial semisimple super-form considered in dimensions $4$ and $8$.

It is immediate that $\LL(\Gamma,\DD)$ is purely even if $|\Gamma|$ is odd.
For $\DD\in\CD(\Gamma)$, the group $\G(\LL(\Gamma,\DD))$ is isomorphic to $\G(\kk^\Gamma)$.
An easy computation shows that the element $\frac12\bigl(v(1,\bar0)+v(1,\bar1)\bigr)\in\LL(\Gamma,\DD)$ is a normalized two-sided integral, and hence we get the following.
\begin{proposition}
$\LL(\Gamma,\DD)$ is semisimple.
\end{proposition}

\begin{definition}
Let $\DD=(\ttt,\ppp;\lambdaa,\muu),\DD'=(\ttt',\ppp';\lambdaa',\muu')$ be two compatible data for $\Gamma$.
We say that $\DD$ is {\it equivalent} to $\DD'$ if there exist a group automorphism $\psi$ of $\Gamma$ and $f\in C^1(\Gamma)$ such that
$\ttt'=\psi\circ\ttt\circ\psi^{-1}$, $\ppp'=\ppp\circ\psi^{-1}$,
\[
\begin{gathered}
\lambdaa'(x) = \lambdaa(\psi^{-1}(x))\,f(\psi^{-1}(x))^{-1}\,f(\psi^{-1}(\ttt'(x)))^{-1}, \text{ and} \\
\muu'(x,y) = \muu(\psi^{-1}(x),\psi^{-1}(y)) \; \partial_1 f(\psi^{-1}(x),\psi^{-1}(y))
\end{gathered}
\]
for $x,y\in\Gamma$.
In this case, we write $(\psi,f) : \DD \xrightarrow{\approx}\DD'$ or simply $\DD\approx \DD'$.
\end{definition}

The relation $\approx$ is an equivalence relation on $\CD(\Gamma)$.
If $\DD\approx\DD'$ in $\CD(\Gamma)$, then
\[
\LL(\Gamma,\DD) \to \LL(\Gamma,\DD');\quad
v(x,\bar{\epsilon}) \mapsto f(x)^\epsilon\, v'(\psi(x),\bar{\epsilon})
\]
gives an isomorphism of Hopf superalgebras,
where
\[
\LL(\Gamma,\DD')=
\bigoplus_{\substack{x\in\Gamma\\ \epsilon\in\{0,1\}}}
\kk v'(x,\bar\epsilon).
\]
Thus we have the following.
\begin{proposition} \label{prp:D=D'}
If $\DD\approx\DD'$ in $\CD(\Gamma)$, then
$\LL(\Gamma,\DD)\cong\LL(\Gamma,\DD')$ as Hopf superalgebras.
\end{proposition}

Let $\DD=(\ttt,\ppp;\lambdaa,\muu)\in\CD(\Gamma)$.
Suppose $[\muu]=[\muu']$ in the second cohomology group $H^2(\Gamma)$.
By definition, there exists $f\in C^1(\Gamma)$ such that $\muu'=\muu\,\partial_1f$.
In this case, if we put $\lambdaa':\Gamma\to\kk^\times; x\mapsto \lambdaa(x)f(x)^{-1}f(\ttt(x))^{-1}$,
then $\DD':=(\ttt,\ppp;\lambdaa',\muu')$ belongs to $\CD(\Gamma)$ and is equivalent to $\DD$.
Hence $\LL(\Gamma,\DD) \cong \LL(\Gamma,\DD')$ by the preceding proposition.

\section{Classification of semisimple non-pointed Hopf superalgebras} \label{sec:class}
Let $\ii$ denote a fixed primitive fourth root of unity in $\kk$.
The aim of this section is to classify semisimple non-pointed Hopf superalgebras of dimension at most $10$ with nonzero odd part.

\subsection{Dimensions \texorpdfstring{$p^2$ and $2p$}{p2 and 2p}}
Let $p$ be an odd prime.
In accordance with the convention fixed in Section~\ref{sec:intro}, write
\[
D_{2p}=\langle a,b\mid a^2=b^p=1,\ aba=b^{-1}\rangle\cong C_p\rtimes C_2
\]
for the dihedral group of order $2p$.

\begin{theorem} \label{thm:main1}
Every Hopf superalgebra $\HH$ of dimension $p^2$ is purely even.
\end{theorem}
\begin{proof}
Suppose first that $\HH$ is non-semisimple, and put $A:=\widehat{\HH}$.
Then $A$ is non-semisimple of dimension $2p^2$.
By~\cite[Theorem~I]{HN09}, either $A$ or $A^*$ is pointed.
Both are non-semisimple, since semisimplicity is preserved by duality for finite-dimensional Hopf algebras in characteristic zero; see~\cite[Theorem~1.1]{HN09}.
By Proposition~\ref{prp:dual-hat} and the properties of bosonization recalled in Subsection~\ref{subsec:Hopf-super}, either $\HH$ or $\HH^*$ is non-semisimple and pointed.
By~\cite[Theorem~4.7]{SSW24}, this Hopf superalgebra is purely even, and hence so is $\HH$.

Now suppose that $\HH$ is semisimple and that $\HH_{\bar1}\neq0$.
Then its bosonization $\widehat{\HH}$ is a semisimple Hopf algebra of dimension $2p^2$.
By Proposition~\ref{prp:super-forms}\eqref{prp:super-forms:(2)}, it is neither commutative nor cocommutative.
By Natale's classification~\cite[Theorem, p.~243]{Nat99}, $\widehat{\HH}$ is isomorphic to the abelian extension
$\mathscr{B}_0^*:=\kk^{D_{2p}}\#_{\sigma,\tau}\kk C_p$ or to its dual $\mathscr{B}_0$.
By~\cite[Lemma~1.3.11]{Nat99}, their groups of group-like elements are
\[
\G(\mathscr{B}_0^*)\cong C_{2p},\qquad
\G(\mathscr{B}_0)\cong C_p^2.
\]
Thus, for either choice, one of $\G(\widehat{\HH})$ and $\G(\widehat{\HH}^*)$ has no direct factor $C_2$, contradicting Proposition~\ref{prp:super-forms}\eqref{prp:super-forms:(1)}.
\end{proof}

By~\cite{Mas96,Ng02}, the purely even classes in dimension $p^2$ are the two group algebras and the $p-1$ Taft algebras.

We next classify the non-trivial semisimple super-forms of dimension $2p$.

By~\cite[Theorems~1.5 and~1.6]{ENO11}, every semisimple Hopf algebra of dimension $4p$ and its dual are of Frobenius type, that is, the dimension of every simple module divides $4p$.
Thus Natale's theorem gives exactly two noncommutative and noncocommutative classes, up to isomorphism, and both are self-dual~\cite[Theorem~5.4.2]{Nat04}.
Following Natale's notation~\cite[Introduction]{Nat01}, write $\mathcal{A}_0$ for the class whose group of group-like elements is $C_2^2$, and $\mathcal{A}_1$ for the class whose group of group-like elements is $C_4$.
Proposition~\ref{prp:super-forms}\eqref{prp:super-forms:(1)} excludes $\mathcal{A}_1$, since $C_4$ is not isomorphic to $C_2\times K$ for any group $K$.
The remaining class $\mathcal{A}_0$ is the abelian extension $\kk^{D_{2p}}\#_{\sigma,\tau}\kk C_2$ used below.

\begin{lemma}
The set $\SD(\mathcal{A}_0)\mathord{/}\mathord{\sim}$ is a singleton.
\end{lemma}
\begin{proof}
We denote the generator of $C_2$ by $t$ as before.
By definition of $\mathcal{A}_0$, the action $\triangleright$, the cocycle $\sigma$, and the dual cocycle $\tau$ are all trivial, and
\[
a\triangleleft t=a,\qquad b\triangleleft t=b^{-1}.
\]
Set
\[
\xi:=\sum_{i=0}^1\sum_{j=0}^{p-1}(-1)^i e_{a^ib^j}\in\kk^{D_{2p}}.
\]
Then $\G(\mathcal{A}_0)=\{\unit,t,\xi,\xi t\}$, while
\[
\mathcal{G}^*(C_2,D_{2p})
 =\{(1,\beta),(a,\beta)\mid \beta(t)\in\{1,-1\}\}.
\]
The admissible data are listed in Table~\ref{tab:A0-admissible-data}.

\begin{table}[H]
\centering
\caption{Admissible data for the abelian extension $\mathcal{A}_0$}
\label{tab:A0-admissible-data}
\begin{tabular}{@{}ccc@{}}
\toprule
$g$ & $\sss$ & $\beta(t)$\\
\midrule
$t$      & $1$ & $-1$\\
$t$      & $a$ & $-1$\\
$\xi$    & $a$ & $1$\\
$\xi$    & $a$ & $-1$\\
$\xi t$  & $1$ & $-1$\\
$\xi t$  & $a$ & $1$\\
\bottomrule
\end{tabular}
\end{table}

Condition (i) of Lemma~\ref{prp:SD(abel)} forces $\ell=1$ and $\sss=a$: conjugation by $a$ is non-trivial, and $b\triangleleft t=b^{-1}\neq b$ since $p$ is odd.
Among the data in Table~\ref{tab:A0-admissible-data}, the super-data are
\[
(t,\alpha_{a,\beta})
\quad\text{and}\quad
(\xi t,\alpha_{a,\beta'}),
\quad\text{where }\beta(t)=-1\text{ and }\beta'(t)=1.
\]
The map
\[
\mathcal{A}_0\longrightarrow\mathcal{A}_0,
\qquad
e_{a^ib^j}t^\ell\longmapsto(-1)^{i\ell}e_{a^ib^j}t^\ell
\]
is a Hopf algebra automorphism and interchanges these two super-data under $\sim$.
Hence $\SD(\mathcal{A}_0)\mathord{/}\mathord{\sim}$ is a singleton.
\end{proof}

Thus, $\mathcal{A}_0$ has a unique non-trivial super-form up to isomorphism.
For the representative $(t,\alpha_{a,\beta(t)=-1})$,
Corollary~\ref{cor:ell-one-identification} identifies this super-form with
$((\kk D_{2p})^a)^*$.
Apply Theorem~\ref{thm:AEG} to the reflection $a\in D_{2p}$.
The resulting Hopf superalgebra $(\kk D_{2p})^a$ is semisimple of dimension $2p$.
Together with the preceding lemma, this gives the following classification.
\begin{theorem} \label{thm:main2}
If $\HH$ is a semisimple non-pointed Hopf superalgebra of dimension $2p$ with $\HH_{\bar1}\neq0$, then
\[
\HH\cong(\kk D_{2p})^a.
\]
Moreover, $\HH$ is self-dual.
\end{theorem}
\begin{proof}
The classification of ordinary semisimple Hopf algebras of dimension $4p$ and the exclusion of $\mathcal{A}_1$ above show that $\widehat{\HH}\cong\mathcal{A}_0$.
The preceding lemma therefore gives at most one isomorphism class of such $\HH$.
For existence, the semisimple Hopf superalgebra $(\kk D_{2p})^a$ has nonzero odd part, since
\[
a(b-b^{-1})a=-(b-b^{-1})\neq0.
\]
It is non-pointed because a semisimple pointed Hopf superalgebra is purely even.
This proves existence and uniqueness.
Since $\mathcal{A}_0$ is self-dual, Proposition~\ref{prp:dual-hat} gives
$\widehat{\HH^*}\cong\widehat{\HH}^{\,*}\cong\mathcal{A}_0$.
The odd part of $\HH^*$ is nonzero, so uniqueness of the non-trivial super-form of $\mathcal{A}_0$ implies $\HH^*\cong\HH$.
\end{proof}

For $p=3$, $D_6\cong\mathfrak{S}_3$, and $\mathcal{A}_0$ is isomorphic to the Hopf algebra $A_+$ of~\cite{Fuk97}.
Thus, the six-dimensional class $(\kk\mathfrak{S}_3)^a$ is isomorphic to
the Hopf superalgebra $\mathcal A_6$ of~\cite[Proposition~5.11 and Theorem~5.12]{SW25}.

\subsection{Dimension \texorpdfstring{$4$}{4}} \label{sec:dim=4}
In~\cite[Section~5.3]{SW25}, the two semisimple non-pointed Hopf superalgebras of dimension $4$ were obtained by computing the coinvariants associated with the super-data of the Kac--Paljutkin Hopf algebra.
Here we treat the same example by the compatible-data construction of Section~\ref{sec:abel}.
We use the same method in dimension $8$.

Let $G = C_2^2 = \langle x,y \mid x^2=y^2=1, xy=yx \rangle$ and let $F=C_2=\langle t\rangle$.
Define
\[
\begin{gathered}
(x^iy^j)\triangleleft t := x^jy^i,\;
\sigma(x^iy^j;t,t) := (-1)^{ij}, \\
(x^iy^j)\triangleright t := t \;(\text{i.e., trivial}),\;
\tau(x^iy^j,x^ry^s;t) := (-1)^{jr}
\end{gathered}
\]
for $i,j,r,s\in\{0,1\}$.
Then we have the abelian extension $\AKP:=\kk^G\#_{\sigma,\tau}\kk F$ of $\kk F$ by $\kk^G$
which is the eight-dimensional Kac--Paljutkin Hopf algebra; see \cite{Mas95b}.

The following is easy to see.
\[
\begin{gathered}
\mathcal{G}(F,G) = \{ (\gamma,1) \mid (\gamma(x),\gamma(y))=(1,1),(-1,1),(1,-1),(-1,-1) \}, \\
\mathcal{G}^*(F,G) = \{ (1,\beta(t)=1),(1,\beta(t)=-1),(xy,\beta(t)=\ii),(xy,\beta(t)=-\ii) \}.
\end{gathered}
\]
To simplify the notation, we fix a Hopf algebra isomorphism
\[
\kappa : \kk^{G} \to \kk G; \quad e_{x^iy^j} \mapsto \frac14 (1+(-1)^ix)(1+(-1)^jy).
\]
Admissibility forces $\sss=xy$ and $\beta(t)^2=-1$.
Equation~\eqref{eq:SD(abel)} then gives $\gamma(x)\gamma(y)=-1$, so $\kappa(g)=x$ or $y$.
This gives the complete list of super-data in Table~\ref{tab:dim4-super-data}.
\begin{table}[H]
\centering
\caption{Super-data for the eight-dimensional abelian extension $\AKP$}
\label{tab:dim4-super-data}
\begin{tabular}{@{}>{$}c<{$}>{$}c<{$}>{$}c<{$}>{$}c<{$}@{}}
\toprule
\kappa(\sum_{z}\gamma(z)e_z) & \fff & \sss & \beta(t)\\
\midrule
y & 1 & xy & -\ii \\
y & 1 & xy & \ii \\
x & 1 & xy & -\ii \\
x & 1 & xy & \ii\\
\bottomrule
\end{tabular}
\end{table}

All four super-data have $\fff=1$ and $G/\langle\sss\rangle\cong C_2$.
Hence Proposition~\ref{prp:induced-compatible-datum} realizes every super-form as $\LL(C_2,\DD)$ for a compatible datum $\DD$.
We first determine the super-data orbits.
The Hopf automorphism group of $\AKP$ was computed in~\cite[Section~4.2, Table~2]{SV12}; see also~\cite{SH25} for automorphism groups of Suzuki's Hopf algebras.
We give its elements in the abelian-extension presentation used here.
\begin{lemma} \label{prp:auto-dim=4}
The cardinality of $\SD(\AKP)\mathord{/}\mathord{\sim}$ is $2$.
\end{lemma}
\begin{proof}
Since the abelian extension $\AKP=\kk^G\#_{\sigma,\tau}\kk F$ satisfies Assumption~\ref{ass:actions},
we can determine the group $\Aut(\AKP)$ of automorphisms of $\AKP$ by Proposition~\ref{prp:automorphism} as follows.
First, determine automorphisms of $G=C_2^2$ which satisfy the condition~\eqref{eq:action-automorphism}.
The only such automorphisms are $\overline{\varphi}^{(1)}$ and $\overline{\varphi}^{(2)}:=\id$,
where $\overline{\varphi}^{(1)}(x)=y, \overline{\varphi}^{(1)}(y)=x$.
Next, solve the equation~\eqref{eq:automorphism}.
For $\overline{\varphi}^{(1)}$ (resp.~$\overline{\varphi}^{(2)}$), the solutions are as follows.
\[
c_{x^iy^j}^{(1,1)}=(-1)^{i+j+ij},\;
c_{x^iy^j}^{(1,2)}=(-1)^{ij} \quad
(\text{resp. }
c_{x^iy^j}^{(2,1)}=(-1)^{i+j},\;
c_{x^iy^j}^{(2,2)}=1),
\]
where $i,j\in\{0,1\}$.
Thus $\Aut(\AKP)$ consists of the following four automorphisms.
\[
\varphi^{(r,s)} : \AKP \to \AKP;\quad
e_{x^iy^j}t^\ell \mapsto (c_{x^iy^j}^{(r,s)})^\ell e_{\overline{\varphi}^{(r)}(x^iy^j)} t^\ell
\qquad(r,s\in\{1,2\}).
\]
Note that $\varphi^{(2,2)}=\id$ and $\Aut(\AKP)=\langle \varphi^{(1,1)},\varphi^{(2,1)}\rangle \cong C_2^2$.
Therefore, by Proposition~\ref{prp:automorphism2},
the set $\{ (\kappa^{-1}(y),\alpha_{xy,\beta(t)=\ii}),\; (\kappa^{-1}(y),\alpha_{xy,\beta(t)=-\ii})\}$
is a complete set of representatives of $\SD(\AKP)\mathord{/}\mathord{\sim}$.
\end{proof}

Identify the generators $x,y,z$ in~\cite[Table~2]{SV12} with $\kappa^{-1}(x)$, $\kappa^{-1}(y)$, and $t$, respectively.
Then their automorphisms $\tau_1$, $\tau_2$, $\tau_3$, and $\tau_4$ correspond to $\varphi^{(2,2)}$, $\varphi^{(2,1)}$, $\varphi^{(1,2)}$, and $\varphi^{(1,1)}$, respectively.

For the two orbit representatives in Lemma~\ref{prp:auto-dim=4}, we have
$\gamma(x)=1$ and $\gamma(y)=-1$. Thus
\[
\Gamma=\operatorname{Ker}(\gamma)=\langle x\rangle\cong C_2.
\]
The construction of Proposition~\ref{prp:induced-compatible-datum} gives
$\ttt=\id$, $\ppp(x)=1$, $\lambdaa(x)=\beta(t)^{-1}$, and
$\muu(x,x)=-1$.
Replacing this datum by the equivalent datum associated with the normalized
cochain $f(x)=\ii$ gives $\muu'=1$ and $\lambdaa'(x)=\beta(t)$.
Hence, on setting $\lambdaa^\pm(x)=\pm\ii$, the two orbit representatives give
\[
\DD_4^\pm=(\id,\ppp;\lambdaa^\pm,1).
\]

\begin{theorem}\label{thm:dim4-classification}
$\LL_4^{(1)}:=\LL(\Gamma,\DD_4^+)$ and $\LL_4^{(2)}:=\LL(\Gamma,\DD_4^-)$ are the two non-trivial semisimple Hopf superalgebras of dimension $4$, up to isomorphism.
They are non-pointed and self-dual.
\end{theorem}
\begin{proof}
By the classification of semisimple Hopf algebras of dimension $8$~\cite{Mas95b} and Proposition~\ref{prp:super-forms}\eqref{prp:super-forms:(2)}, the bosonization of any non-trivial semisimple Hopf superalgebra of dimension $4$ is isomorphic to $\AKP$.
Lemma~\ref{prp:auto-dim=4} and Proposition~\ref{prp:one-to-one} show that the two classes are distinct and exhaust all non-trivial super-forms of $\AKP$.
They are non-pointed because their odd parts are nonzero and every semisimple pointed Hopf superalgebra is purely even.
In the notation of~\cite{SW25}, they are $\mathcal{A}_4(\ii)$ and $\mathcal{A}_4(-\ii)$, respectively.
Their self-duality follows from~\cite[Theorem~5.10]{SW25}.
\end{proof}

\begin{remark}\label{rem:SW25-self-duality}
For $\HH=\AA_4(-\ii)$, let $1^*,x^*,(x^2)^*,z^*$ be the basis of $\HH^*$ dual to $1,x,x^2,z$.
Choose $c\in\kk$ with $c^2=\ii$.
By~\cite[Section~5.3.3]{SW25}, the map
\[
 \HH\longrightarrow\HH^*;\qquad x\longmapsto1^*,\quad z\longmapsto cz^*
\]
is an isomorphism of Hopf superalgebras.
Replacing $\ii$ by $-\ii$ gives the other case.
\end{remark}

\subsection{Dimension \texorpdfstring{$8$}{8}} \label{sec:dim=8}
Kashina classified the semisimple Hopf algebras of dimension $16$ that are neither commutative nor cocommutative~\cite[Theorem~1.2]{Kas00}.
By Proposition~\ref{prp:super-forms}\eqref{prp:super-forms:(2)}, the bosonization of a non-trivial semisimple Hopf superalgebra of dimension $8$ belongs to this list.
Table~\ref{tab:kashina16} gives these sixteen classes and their groups of group-like elements~\cite[Table~1]{Kas00}.
The dual identifications in the third and fourth rows are justified in Section~\ref{sec:dual}.
\begin{table}[H]
\centering
\caption{Noncommutative and noncocommutative semisimple Hopf algebras of dimension $16$ in Kashina's classification}
\label{tab:kashina16}
\begin{tabular}{@{}lcc@{}}
\toprule
$A$ & $\G(A)$ & $\G(A^*)$\\
\midrule
$H_{d:-1,1}\cong\AKP\otimes\kk C_2$ & $C_2^3$ & $C_2^3$\\
$H_{d:1,1}\cong\kk(D_8\times C_2)_J$ & $C_2^3$ & $C_2^3$\\
$H_{d:-1,-1}\cong (H_{c:\sigma_1})^*$ & $C_2^3$ & $C_2\times C_4$\\
$H_{d:1,-1}\cong (H_{b:1})^*$ & $C_2^3$ & $C_2\times C_4$\\
\addlinespace
$H_{c:\sigma_1}$ & $C_2\times C_4$ & $C_2^3$\\
$H_{b:1}$ & $C_2\times C_4$ & $C_2^3$\\
$H_{c:\sigma_0}$ & $C_2\times C_4$ & $C_2\times C_4$\\
$H_{a:1}$ & $C_2\times C_4$ & $C_2\times C_4$\\
$H_{a:y}$ & $C_2\times C_4$ & $C_2\times C_4$\\
$H_{b:y}$ & $C_2\times C_4$ & $C_2\times C_4$\\
$H_{b:x^2y}$ & $C_2\times C_4$ & $C_2\times C_4$\\
\addlinespace
$H_{C:1}\cong(\kk D_{16})_J$ & $D_8$ & $C_2^2$\\
$H_E\cong(\kk G_2)_J$ & $D_8$ & $C_2^2$\\
$H_{B:1}\cong((\kk D_{16})_J)^*$ & $C_2^2$ & $D_8$\\
$H_{B:X}\cong((\kk G_2)_J)^*$ & $C_2^2$ & $D_8$\\
$H_{C:\sigma_1}$ & $C_2^2$ & $C_2^2$\\
\bottomrule
\end{tabular}
\end{table}
We follow Kashina's notation throughout this subsection.
Here, $\AKP$ is the semisimple non-pointed Hopf algebra of dimension $8$ discussed in Subsection~\ref{sec:dim=4}.

By Proposition~\ref{prp:super-forms}\eqref{prp:super-forms:(1)}, only the entries for which both $\G(A)$ and $\G(A^*)$ have a direct factor isomorphic to $C_2$ can be bosonizations.
This excludes $H_{C:1}$, $H_E$, $H_{B:1}$, and $H_{B:X}$: each has $D_8$ as one of these two groups, and $D_8$ has no direct factor isomorphic to $C_2$.
We examine the remaining twelve entries below.
They are realized by abelian extensions $A=\kk^G\#_{\sigma,\tau}\kk F$ of $\kk F$ by $\kk^G$, where $F=C_2=\langle t\rangle$ and the left action $\triangleright:G\times F\to F$ is trivial.

Since $F=\langle t\rangle$,
for $(\sss,\beta)\in\mathcal{G}^*(F,G)$,
writing $\alpha_{\sss,\beta(t)}$ in place of $\alpha_{\sss,\beta}$ causes no confusion.
Therefore, we use this simplified notation from now on.

\subsubsection{Determination of super-data}
Some of the Hopf algebras in Table~\ref{tab:kashina16} have also been studied as coradicals of larger finite-dimensional Hopf algebras.
Examples include $H_{b:1}$, $H_{d:-1,1}$, $(H_{b:1})^*$, $H_{b:x^2y}$, $H_{a:y}$, and $H_{b:y}$; see~\cite{ZGH21a,ZGH21b,ZGHS23,WWH26,ZY25,Zhe26}, respectively.
Those works study Nichols algebras and liftings, whereas here we determine super-forms.
We give the required automorphisms in the abelian-extension presentations below and use them to calculate the orbits of super-data.

\paragraph{$\G(A)= C_4\times C_2$ case}
We treat the seven Hopf algebras
\[
A\in\{H_{a:1},H_{a:y},H_{b:1},H_{b:y},H_{b:x^2y},
H_{c:\sigma_0},H_{c:\sigma_1}\}.
\]
They are realized as $A=\kk^G\#_{\sigma,\tau}\kk F$, where
\[
G=C_4\times C_2=\langle x,y\mid x^4=y^2=1,\ xy=yx\rangle,
\qquad F=\langle t\rangle.
\]
For these algebras, the dual $2$-cocycle $\tau:G\times G\to (\kk^F)^\times$ is given by
\[
\tau(x^iy^j, x^ry^s;t) = (-1)^{jr} \quad (i,r\in\{0,1,2,3\},\; j,s\in\{0,1\}).
\]

\begin{lemma}\label{lem:five-empty-superdata}
For $A\in\{H_{a:1},H_{a:y},H_{b:1},H_{b:y},H_{b:x^2y}\}$, we have $\SD(A)=\varnothing$.
Moreover, $\AD(H_{b:x^2y})=\varnothing$.
\end{lemma}
\begin{proof}
Suppose $A=H_{a:1}$.
Then the right action is $x^ry^s \triangleleft t = x^{r+2s}y^s$ and the $2$-cocycle is $\sigma(x^ry^s;t,t) = 1$ ($r\in\{0,1,2,3\},s\in\{0,1\}$).
Suppose $A=H_{a:y}$.
Then the right action is $x^ry^s \triangleleft t = x^{r+2s}y^s$ and the $2$-cocycle is $\sigma(x^ry^s;t,t) = (-1)^s$ ($r\in\{0,1,2,3\},s\in\{0,1\}$).

Suppose $A=H_{b:1}$.
Then the right action is $x^ry^s \triangleleft t = x^{-r}y^s$ and the $2$-cocycle is $\sigma(x^ry^s;t,t) = 1$ ($r\in\{0,1,2,3\},s\in\{0,1\}$).
Suppose $A=H_{b:y}$.
Then the right action is $x^ry^s \triangleleft t = x^{-r}y^s$ and the $2$-cocycle is $\sigma(x^ry^s;t,t) = (-1)^s$ ($r\in\{0,1,2,3\},s\in\{0,1\}$).
Suppose $A=H_{b:x^2y}$.
Then the right action is $x^ry^s \triangleleft t = x^{-r}y^s$ and the $2$-cocycle is $\sigma(x^ry^s;t,t) = (-1)^{r+s}$ ($r\in\{0,1,2,3\},s\in\{0,1\}$).

For all five algebras above, every group-like element lies in $\kk^G$, since $\tau(x,y;t)=1\neq-1=\tau(y,x;t)$.
If $g=\sum_{u\in G}\gamma(u)e_u$ has order $2$, then $\gamma(x),\gamma(y)\in\{1,-1\}$, and the formulas for the actions give $\gamma(u\triangleleft t)=\gamma(u)$ for every $u\in G$.
Thus every group-like element of order $2$ is central, and $\SD(A)=\varnothing$.
For $H_{b:x^2y}$, admissibility moreover forces $\sss\in\{1,x^2\}$, whereas every character of order dividing $2$ takes value $1$ on these elements; hence $\AD(H_{b:x^2y})=\varnothing$.
\end{proof}

By Lemma~\ref{lem:five-empty-superdata}, it remains to treat $H_{c:\sigma_0}$ and $H_{c:\sigma_1}$ in this case.
Suppose $A=H_{c:\sigma_0}$.
Then the right action is $x^ry^s \triangleleft t = x^{r}y^{r+s}$ and the $2$-cocycle is $\sigma(x^ry^s;t,t) = (-1)^{r(r-1)/2}$ ($r\in\{0,1,2,3\},s\in\{0,1\}$).
If we fix a Hopf algebra isomorphism
\[
\kappa : \kk^{G} \to \kk G; \quad e_{x^ry^s} \mapsto \frac18 (1+\ii^r x + (-1)^rx^2 + (-\ii)^rx^3)(1+(-1)^sy),
\]
then the super-data are precisely those listed in Table~\ref{tab:SD-Hc-sigma0}.
\begin{table}[H]
\centering
\caption{Super-data for $H_{c:\sigma_0}$}
\label{tab:SD-Hc-sigma0}
\begin{tabular}{@{}>{$}c<{$}>{$}c<{$}>{$}c<{$}>{$}c<{$}@{}}
\toprule
\kappa(\sum_{z}\gamma(z)e_z) & \fff & \sss & \beta(t)\\
\midrule
x^2y & 1 & y & -1 \\
x^2y & 1 & y & 1 \\
y & 1 & y & -1 \\
y & 1 & y & 1\\
\bottomrule
\end{tabular}
\end{table}
\begin{lemma} \label{prp:Aut(H_c:sigma_0)}
For $A=H_{c:\sigma_0}$,
\[
\SD(A)\mathord{/}\mathord{\sim}
 = \{ [(\kappa^{-1}(x^2y),\alpha_{y,-1})] \}.
\]
\end{lemma}
\begin{proof}
By Lemma~\ref{lem:automorphism-count}, the induced group automorphism $\overline\varphi$ commutes with $u\mapsto u\triangleleft t$.
Since $x\triangleleft t=xy$ and $y\triangleleft t=y$, this gives
\[
\overline\varphi(x)\in\{x,x^{-1},xy,x^{-1}y\},
\qquad \overline\varphi(y)=y.
\]
For each of these four possibilities, any two solutions of~\eqref{eq:automorphism} differ by a character $\chi$ satisfying $\chi(x)^2\chi(y)=1$.
There are four such characters: $\chi(x)$ is an arbitrary fourth root of unity and $\chi(y)=\chi(x)^2$.
Consequently, $|\Aut(A)|\le4\cdot4=16$.
The following maps satisfy equations~\eqref{eq:action-automorphism} and~\eqref{eq:automorphism}, and hence are Hopf algebra automorphisms, where $r\in\{0,1,2,3\}$ and $s,\ell\in\{0,1\}$:
\[
\begin{aligned}
\varphi_1(e_{x^ry^s}t^\ell)
 &=(-1)^{s\ell}\ii^{r\ell}e_{x^ry^s}t^\ell,\\
\varphi_2(e_{x^ry^s}t^\ell)
 &=(-1)^{r\ell}\ii^{r\ell}e_{x^{-r}y^s}t^\ell,\\
\varphi_3(e_{x^ry^s}t^\ell)
 &=(-1)^{(s+r(r-1)/2)\ell}e_{x^{-r}y^{r+s}}t^\ell.
\end{aligned}
\]
Here and below, exponents of $x$ and $y$ are read modulo $4$ and $2$, respectively.
The first two generators satisfy
$\varphi_1^4=\varphi_2^2=\id$ and
$\varphi_2\varphi_1\varphi_2=\varphi_1^{-1}$,
while $\varphi_3$ is a central involution.
The $16$ products $\varphi_1^i\varphi_2^j\varphi_3^k$ ($0\le i<4$, $0\le j,k<2$) are distinct, so they attain the upper bound. Hence
\[
 \Aut(A)=\langle\varphi_1,\varphi_2\rangle
          \times\langle\varphi_3\rangle
 \cong D_8\times C_2.
\]
Applying Proposition~\ref{prp:automorphism2} to these generators shows that their action on the four super-data in Table~\ref{tab:SD-Hc-sigma0} is transitive.
Thus the super-datum in the statement represents the unique orbit.
\end{proof}

Suppose $A=H_{c:\sigma_1}$.
Then the right action is $x^ry^s \triangleleft t = x^{r}y^{r+s}$ and the $2$-cocycle is $\sigma(x^ry^s;t,t) = (-1)^{r(r-1)/2}\ii^r$ ($r\in\{0,1,2,3\},s\in\{0,1\}$).
The super-data are precisely those listed in Table~\ref{tab:SD-Hc-sigma1}.
\begin{table}[H]
\centering
\caption{Super-data for $H_{c:\sigma_1}$}
\label{tab:SD-Hc-sigma1}
\begin{tabular}{@{}>{$}c<{$}>{$}c<{$}>{$}c<{$}>{$}c<{$}@{}}
\toprule
\kappa(\sum_{z}\gamma(z)e_z) & \fff & \sss & \beta(t)\\
\midrule
x^2y & 1 & y & -1 \\
x^2y & 1 & y & 1 \\
x^2y & 1 & x^2y & -1 \\
x^2y & 1 & x^2y & 1 \\
y & 1 & y & -1 \\
y & 1 & y & 1 \\
y & 1 & x^2y & -1 \\
y & 1 & x^2y & 1\\
\bottomrule
\end{tabular}
\end{table}
\begin{lemma} \label{prp:Aut(H_c:sigma_1)}
For $A=H_{c:\sigma_1}$, the following three elements form a complete set of representatives of $\SD(A)\mathord{/}\mathord{\sim}$:
\[
(\kappa^{-1}(x^2y),\alpha_{y,-1}),\qquad
(\kappa^{-1}(x^2y),\alpha_{x^2y,-1}),\qquad
(\kappa^{-1}(x^2y),\alpha_{x^2y,1}).
\]
\end{lemma}
\begin{proof}
By Lemma~\ref{lem:automorphism-count}, the induced group automorphism $\overline\varphi$ commutes with $u\mapsto u\triangleleft t$.
Since $x\triangleleft t=xy$ and $y\triangleleft t=y$, this gives
\[
\overline\varphi(x)\in\{x,x^{-1},xy,x^{-1}y\},
\qquad \overline\varphi(y)=y.
\]
For each of these four possibilities, any two solutions of~\eqref{eq:automorphism} differ by a character $\chi$ satisfying $\chi(x)^2\chi(y)=1$.
There are four such characters: $\chi(x)$ is an arbitrary fourth root of unity and $\chi(y)=\chi(x)^2$.
Consequently, $|\Aut(A)|\le4\cdot4=16$.
The following maps satisfy equations~\eqref{eq:action-automorphism} and~\eqref{eq:automorphism}:
\[
\begin{aligned}
\psi_1(e_{x^ry^s}t^\ell)
 &=(-1)^{s\ell}\ii^{r\ell}e_{x^ry^s}t^\ell,\\
\psi_2(e_{x^ry^s}t^\ell)
 &=e_{x^{-r}y^s}t^\ell,\\
\psi_3(e_{x^ry^s}t^\ell)
 &=(-1)^{s\ell}\ii^{r^2\ell}e_{x^{-r}y^{r+s}}t^\ell,
\end{aligned}
\]
where $r\in\{0,1,2,3\}$ and $s,\ell\in\{0,1\}$.
The relations are
$\psi_1^4=\psi_2^2=\psi_3^2=\id$,
$\psi_2\psi_1\psi_2=\psi_1^{-1}$, and
$\psi_3\psi_j=\psi_j\psi_3$ for $j=1,2$.
The $16$ products $\psi_1^i\psi_2^j\psi_3^k$ ($0\le i<4$, $0\le j,k<2$) are distinct and attain the upper bound.
Consequently,
\[
 \Aut(A)=\langle\psi_1,\psi_2\rangle
          \times\langle\psi_3\rangle
 \cong D_8\times C_2.
\]
Number the super-data in Table~\ref{tab:SD-Hc-sigma1} from $1$ to $8$. By Proposition~\ref{prp:automorphism2}, the orbits are $\{1,2,5,6\}$, $\{3,8\}$, and $\{4,7\}$.
The three elements in the statement represent these orbits in the same order.
\end{proof}

\paragraph{$\G(A)= C_2^3$ case}
We treat
\[
A\in\{H_{d:1,1},H_{d:-1,1},H_{d:1,-1},H_{d:-1,-1}\},
\]
realized as $A=\kk^G\#_{\sigma,\tau}\kk F$, where
\[
G=C_2^3=\langle y,z,w\mid y^2=z^2=w^2=1,\ yz=zy,\ yw=wy,\ zw=wz\rangle,
\qquad F=\langle t\rangle.
\]
In this case, the right action is given by $y^qz^sw^u \triangleleft t = y^sz^qw^u$ ($q,s,u\in\{0,1\}$).
For $\epsilon,\eta\in\{\pm1\}$,
the $2$-cocycle $\sigma_{\epsilon,\eta} : F\times F\to (\kk^G)^\times$ and the dual $2$-cocycle $\tau_{\epsilon,\eta} : G\times G\to (\kk^F)^\times$
of $H_{d:\epsilon,\eta}$ are given as
\[
\begin{gathered}
\sigma(t,t) = e_{1} + e_{y} + e_{z} + \epsilon e_{yz} + \eta (e_w + e_{yw} + e_{zw}) + \epsilon\eta e_{yzw}, \\
\tau(y^qz^sw^u,y^{q'}z^{s'}w^{u'};t) = \epsilon^{sq'}(-1)^{u(q'+s')} \quad (q,s,u,q',s',u'\in\{0,1\}),
\end{gathered}
\]
respectively.
We fix a Hopf algebra isomorphism
\[
\kappa : \kk^{G} \to \kk G; \quad e_{y^qz^sw^u} \mapsto \frac18 (1+(-1)^qy)(1+(-1)^sz)(1+(-1)^uw).
\]

Suppose $A=H_{d:1,1}$.
The super-data are precisely those listed in Table~\ref{tab:SD-Hd-11}.
\begin{table}[H]
\centering
\caption{Super-data for $H_{d:1,1}$}
\label{tab:SD-Hd-11}
\begin{tabular}{@{}>{$}c<{$}>{$}c<{$}>{$}c<{$}>{$}c<{$}@{}}
\toprule
\kappa(\sum_{z}\gamma(z)e_z) & \fff & \sss & \beta(t)\\
\midrule
zw & 1 & w & -1 \\
zw & 1 & w & 1 \\
yw & 1 & w & -1 \\
yw & 1 & w & 1 \\
z & 1 & yzw & -1 \\
z & 1 & yzw & 1 \\
y & 1 & yzw & -1 \\
y & 1 & yzw & 1\\
\bottomrule
\end{tabular}
\end{table}
\begin{lemma} \label{prp:Aut(H_d:1,1)}
For $A=H_{d:1,1}$,
\[
\SD(A)\mathord{/}\mathord{\sim}
 =\{ [(\kappa^{-1}(zw), \alpha_{w,-1})] \}.
\]
\end{lemma}
\begin{proof}
By Lemma~\ref{lem:automorphism-count}, the induced group automorphism $\overline\varphi$ commutes with the action interchanging $y$ and $z$ and fixing $w$.
The group automorphisms with this property are given by
\[
\overline\varphi(y)\in\{y,z,yw,zw\},\qquad
\overline\varphi(z)=\overline\varphi(y)\triangleleft t,\qquad
\overline\varphi(w)\in\{w,yzw\},
\]
so there are eight possibilities.
For each, any two solutions of~\eqref{eq:automorphism} differ by a character $\chi$ with $\chi(y)=\chi(z)$ and arbitrary $\chi(w)\in\{\pm1\}$.
There are four such characters, so $|\Aut(A)|\le8\cdot4=32$.
For $q,s,u,\ell\in\{0,1\}$, the following maps satisfy equations~\eqref{eq:action-automorphism} and~\eqref{eq:automorphism}:
\[
\begin{aligned}
\varphi_1(e_{y^qz^sw^u}t^\ell)
 &=(-1)^{(q+s)\ell}e_{y^qz^sw^u}t^\ell,\\
\varphi_2(e_{y^qz^sw^u}t^\ell)
 &=(-1)^{u\ell}e_{y^sz^qw^u}t^\ell,\\
\varphi_3(e_{y^qz^sw^u}t^\ell)
 &=e_{y^sz^qw^u}t^\ell,\\
\varphi_4(e_{y^qz^sw^u}t^\ell)
 &=e_{y^{s+u}z^{q+u}w^u}t^\ell,\\
\varphi_5(e_{y^qz^sw^u}t^\ell)
 &=(-1)^{qs\ell}\ii^{(q-s)\ell}
   e_{y^sz^qw^{q+s+u}}t^\ell.
\end{aligned}
\]
All exponents of $y,z,w$ are read modulo $2$.
The four automorphisms $\varphi_1,\ldots,\varphi_4$ are commuting involutions and generate a normal subgroup
\[
 N:=\langle\varphi_1,\varphi_2,\varphi_3,\varphi_4\rangle\cong C_2^4.
\]
Moreover, $\varphi_5^2=\id$, and conjugation by $\varphi_5$ is determined by
\[
\varphi_5\varphi_1\varphi_5^{-1}=\varphi_1,
\qquad
\varphi_5\varphi_2\varphi_5^{-1}=\varphi_2,
\qquad
\varphi_5\varphi_3\varphi_5^{-1}=\varphi_1\varphi_3,
\qquad
\varphi_5\varphi_4\varphi_5^{-1}=\varphi_2\varphi_4.
\]
Every element of $N$ preserves the exponent $u$ of $w$, whereas $\varphi_5$ changes it to $q+s+u$.
Hence $\varphi_5\notin N$, so these generators give $32$ distinct automorphisms and attain the upper bound. Thus
\[
 \Aut(A)=N\rtimes\langle\varphi_5\rangle
 \cong C_2^4\rtimes C_2,
\]
where the action of the last factor is specified by the preceding conjugation formulas.
Proposition~\ref{prp:automorphism2} now shows that this group acts transitively on the eight super-data in Table~\ref{tab:SD-Hd-11}.
Hence the super-datum in the statement represents the unique orbit.
\end{proof}

For $A=H_{d:1,-1}$, every group-like element lies in $\kk^G$, since $\tau(y,w;t)=1\neq-1=\tau(w,y;t)$.
A super-datum must have $\sss\in\{1,yz,w,yzw\}$, and Lemma~\ref{prp:SD(abel)} excludes $\sss=1,yz$: either choice would imply that the group-like element is central.
Thus $\sss=w$ or $yzw$, but then admissibility requires
\[
\beta(t)^2=\sigma(\sss;t,t)=-1,
\qquad
\beta(t)^2=\tau(\sss,\sss;t)^{-1}=1,
\]
a contradiction. Therefore $\SD(H_{d:1,-1})=\varnothing$.
Suppose $A=H_{d:-1,1}$.
The super-data are precisely those listed in Table~\ref{tab:SD-Hd-m11}.
\begin{table}[H]
\centering
\caption{Super-data for $H_{d:-1,1}$}
\label{tab:SD-Hd-m11}
\begin{tabular}{@{}>{$}c<{$}>{$}c<{$}>{$}c<{$}>{$}c<{$}@{}}
\toprule
\kappa(\sum_{z}\gamma(z)e_z) & \fff & \sss & \beta(t)\\
\midrule
zw & 1 & w & -1 \\
zw & 1 & w & 1 \\
zw & 1 & yz & -\ii \\
zw & 1 & yz & \ii \\
yw & 1 & w & -1 \\
yw & 1 & w & 1 \\
yw & 1 & yz & -\ii \\
yw & 1 & yz & \ii \\
z & 1 & yz & -\ii \\
z & 1 & yz & \ii \\
y & 1 & yz & -\ii \\
y & 1 & yz & \ii\\
\bottomrule
\end{tabular}
\end{table}
\begin{lemma} \label{prp:Aut(H_d:-1,1)}
For $A=H_{d:-1,1}$, a complete set of representatives of
$\SD(A)\mathord{/}\mathord{\sim}$ is
\[
(\kappa^{-1}(zw),\alpha_{w,-1}),\quad
(\kappa^{-1}(zw),\alpha_{yz,-\ii}),\quad
(\kappa^{-1}(zw),\alpha_{yz,\ii}),\quad
(\kappa^{-1}(z),\alpha_{yz,-\ii}).
\]
\end{lemma}
\begin{proof}
By Lemma~\ref{lem:automorphism-count}, commutation with the action of $t$ gives
\[
\overline\varphi(y)\in\{y,z,yw,zw\},\qquad
\overline\varphi(z)=\overline\varphi(y)\triangleleft t,\qquad
\overline\varphi(w)\in\{w,yzw\}.
\]
Comparing the last equation of~\eqref{eq:automorphism} for $(y,w)$ and $(w,y)$ rules out $\overline\varphi(w)=yzw$: the ratio $\tau(y,w;t)/\tau(w,y;t)$ is $-1$, whereas its value at $(\overline\varphi(y),yzw)$ is $1$.
Thus there are at most four possibilities for $\overline\varphi$.
For each, any two solutions of~\eqref{eq:automorphism} differ by a character $\chi$ with $\chi(y)=\chi(z)$ and arbitrary $\chi(w)\in\{\pm1\}$.
There are four such characters, so $|\Aut(A)|\le4\cdot4=16$.
The following maps satisfy equations~\eqref{eq:action-automorphism} and~\eqref{eq:automorphism}:
\[
\begin{aligned}
\varphi_1(e_{y^qz^sw^u}t^\ell)
 &=(-1)^{(s+qs)\ell}\ii^{(q+s)\ell}
   e_{y^qz^sw^{q+s+u}}t^\ell,\\
\varphi_2(e_{y^qz^sw^u}t^\ell)
 &=(-1)^{u\ell}e_{y^qz^sw^u}t^\ell,\\
\varphi_3(e_{y^qz^sw^u}t^\ell)
 &=(-1)^{(u+qs)\ell}e_{y^sz^qw^u}t^\ell,
\end{aligned}
\]
where $q,s,u,\ell\in\{0,1\}$.
Here $\varphi_1$ has order $4$, while $\varphi_2$ and $\varphi_3$ have order $2$;
moreover, $\varphi_2\varphi_1\varphi_2=\varphi_1^{-1}$, and $\varphi_3$ commutes with both $\varphi_1$ and $\varphi_2$.
The $16$ products $\varphi_1^i\varphi_2^j\varphi_3^k$ ($0\le i<4$, $0\le j,k<2$) are distinct and attain the upper bound.
Therefore,
\[
 \Aut(A)=\langle\varphi_1,\varphi_2\rangle
          \times\langle\varphi_3\rangle
 \cong D_8\times C_2.
\]
Number the super-data in Table~\ref{tab:SD-Hd-m11} from $1$ to $12$. By Proposition~\ref{prp:automorphism2}, the orbits are
\[
\{1,2,5,6\},\qquad \{3,8\},\qquad \{4,7\},\qquad \{9,10,11,12\}.
\]
They correspond to the representatives in the statement in the same order, so their sizes are $4$, $2$, $2$, and $4$.
\end{proof}

Suppose $A=H_{d:-1,-1}$.
The super-data are precisely those listed in Table~\ref{tab:SD-Hd-mm}.
\begin{table}[H]
\centering
\caption{Super-data for $H_{d:-1,-1}$}
\label{tab:SD-Hd-mm}
\begin{tabular}{@{}>{$}c<{$}>{$}c<{$}>{$}c<{$}>{$}c<{$}@{}}
\toprule
\kappa(\sum_{z}\gamma(z)e_z) & \fff & \sss & \beta(t)\\
\midrule
zw & 1 & yz & -\ii \\
zw & 1 & yz & \ii \\
yw & 1 & yz & -\ii \\
yw & 1 & yz & \ii \\
z & 1 & yz & -\ii \\
z & 1 & yz & \ii \\
y & 1 & yz & -\ii \\
y & 1 & yz & \ii\\
\bottomrule
\end{tabular}
\end{table}
\begin{lemma} \label{prp:Aut(H_d:-1,-1)}
For $A=H_{d:-1,-1}$, the following three elements form a complete set of representatives of $\SD(A)\mathord{/}\mathord{\sim}$:
\[
(\kappa^{-1}(zw),\alpha_{yz,-\ii}),\qquad
(\kappa^{-1}(z),\alpha_{yz,-\ii}),\qquad
(\kappa^{-1}(z),\alpha_{yz,\ii}).
\]
\end{lemma}
\begin{proof}
By Lemma~\ref{lem:automorphism-count}, commutation with the action of $t$ gives
\[
\overline\varphi(y)\in\{y,z,yw,zw\},\qquad
\overline\varphi(z)=\overline\varphi(y)\triangleleft t,\qquad
\overline\varphi(w)\in\{w,yzw\}.
\]
Comparing the last equation of~\eqref{eq:automorphism} for $(y,w)$ and $(w,y)$ rules out $\overline\varphi(w)=yzw$: the ratio $\tau(y,w;t)/\tau(w,y;t)$ is $-1$, whereas its value at $(\overline\varphi(y),yzw)$ is $1$.
Thus there are at most four possibilities for $\overline\varphi$.
For each, any two solutions of~\eqref{eq:automorphism} differ by a character $\chi$ with $\chi(y)=\chi(z)$ and arbitrary $\chi(w)\in\{\pm1\}$.
There are four such characters, so $|\Aut(A)|\le4\cdot4=16$.
The following maps satisfy equations~\eqref{eq:action-automorphism} and~\eqref{eq:automorphism}:
\[
\begin{aligned}
\varphi_1(e_{y^qz^sw^u}t^\ell)
 &=(-1)^{(q+s+qs)\ell}\ii^{(q+s)\ell}
   e_{y^qz^sw^{q+s+u}}t^\ell,\\
\varphi_2(e_{y^qz^sw^u}t^\ell)
 &=(-1)^{(u+qs)\ell}e_{y^sz^qw^u}t^\ell,\\
\varphi_3(e_{y^qz^sw^u}t^\ell)
 &=(-1)^{(q+s+qs)\ell}e_{y^sz^qw^u}t^\ell,
\end{aligned}
\]
where $q,s,u,\ell\in\{0,1\}$.
The relations are
$\varphi_1^4=\varphi_2^2=\varphi_3^2=\id$,
$\varphi_2\varphi_1\varphi_2=\varphi_1^{-1}$, and
$\varphi_3\varphi_1=\varphi_1\varphi_3$, $\varphi_3\varphi_2=\varphi_2\varphi_3$.
The $16$ products $\varphi_1^i\varphi_2^j\varphi_3^k$ ($0\le i<4$, $0\le j,k<2$) are distinct and attain the upper bound.
Hence
\[
 \Aut(A)=\langle\varphi_1,\varphi_2\rangle
          \times\langle\varphi_3\rangle
 \cong D_8\times C_2.
\]
Number the super-data in Table~\ref{tab:SD-Hd-mm} from $1$ to $8$. By Proposition~\ref{prp:automorphism2}, the orbits are $\{1,2,3,4\}$, $\{5,8\}$, and $\{6,7\}$.
The three elements in the statement represent these orbits in the same order.
\end{proof}

\paragraph{$\G(A)= C_2^2$ case}
We treat $A=H_{C:\sigma_1}=\kk^G\#_{\sigma,\tau}\kk F$, where
\[
G=D_8=\langle a,b\mid a^2=b^4=1,\ aba=b^{-1}\rangle,
\qquad F=C_2=\langle t\rangle.
\]
Using the normal form $b^ra^s$, put
\[
X:=\sum_{r=0}^3\sum_{s=0}^1(-1)^r e_{b^ra^s},
\qquad
Y:=\sum_{r=0}^3\sum_{s=0}^1(-1)^s e_{b^ra^s}.
\]
Then $\G(A)=\{1,X,Y,XY\}\cong C_2^2$.
The dual group-likes are generated by the algebra maps
\[
\chi(e_{b^ra^s}t^\ell)=\delta_{r,2}\delta_{s,0},
\qquad
\varphi(e_{b^ra^s}t^\ell)=(-1)^\ell\delta_{r,0}\delta_{s,0},
\]
where $r\in\{0,1,2,3\}$ and $s,\ell\in\{0,1\}$; hence
$\G(A^*)=\langle\chi\rangle\times\langle\varphi\rangle\cong C_2^2$.
Moreover,
\[
\chi(X)=\chi(Y)=\varphi(X)=\varphi(Y)=1.
\]
Thus every $\alpha\in\G(A^*)$ satisfies $\alpha(g)=1$ for every $g\in\G(A)$.
The admissibility condition $\alpha(g)=-1$ can never hold, so $\AD(A)=\varnothing$ and therefore $\SD(A)=\varnothing$.

\subsubsection{Compatible data attached to the super-data}\label{subsec:compatible-data}
The five Hopf algebras
\[
H_{c:\sigma_0},\ H_{c:\sigma_1},\ H_{d:1,1},\ H_{d:-1,1},\ H_{d:-1,-1}
\]
are precisely the bosonizations for which the preceding calculation gives super-data.
For every super-datum listed above, $\fff=1$, equivalently $\ell=0$ in Section~\ref{sec:ell=0}.
Thus Proposition~\ref{prp:induced-compatible-datum} produces a compatible datum on $C_4$ or $C_2^2$.
We use the following notation for $1$-cochains.

\begin{notation}\label{not:cochain-values}
For a map $f:C_4\to\kk^\times$, where $C_4=\langle x\mid x^4=1\rangle$, write
\[
(f)_{C_4}:=
\bigl(f(1),f(x),f(x^2),f(x^3)\bigr).
\]
For a map $f:C_2^2\to\kk^\times$, where $C_2^2=\langle y,z\mid y^2=z^2=1,\ yz=zy\rangle$, write
\[
(f)_{C_2^2}:=
\bigl(f(1),f(y),f(z),f(yz)\bigr).
\]
\end{notation}

\paragraph{$\Gamma=C_4$ case}
Let $\Gamma=C_4=\langle x\mid x^4=1\rangle$, let $\ppp(x)=1$, and let $\ttt$ denote inversion, $\ttt(x)=x^3$.
We use the tuple notation $(\lambdaa)_{C_4}$ of Notation~\ref{not:cochain-values}.
The four super-data orbits from Lemmas~\ref{prp:Aut(H_c:sigma_0)} and~\ref{prp:Aut(H_c:sigma_1)} yield the following data:
\[
\begin{gathered}
\DD_8^{(1)}=( \id, \ppp; \bigl(1,-1,-1,1\bigr), 1 ), \quad
\DD_8^{(2)}=( \id, \ppp; \bigl(1,-\ii,1,-\ii\bigr), 1 ),\\
\DD_8^{(3)}=( \ttt, \ppp; \bigl(1,\ii,1,\ii\bigr), 1 ), \quad
\DD_8^{(4)}=( \ttt, \ppp; \bigl(1,-\ii,1,-\ii\bigr), 1 ).
\end{gathered}
\]
Set $\LL_8^{(i)}:=\LL(\Gamma,\DD_8^{(i)})$ for $i\in\{1,2,3,4\}$.
Note that $\G(\LL_8^{(i)}) \cong \G(\kk^\Gamma) \cong C_4$ since $\kk^\Gamma\cong \kk \Gamma$.

\begin{proposition} \label{prp:bosonization1}
The four Hopf superalgebras $\LL_8^{(i)}$ ($1\leq i\leq4$) are pairwise non-isomorphic.
Their bosonizations are
\[
\widehat{\LL_8^{(1)}}\cong H_{c:\sigma_0},
\qquad
\widehat{\LL_8^{(i)}}\cong H_{c:\sigma_1}
\quad(i=2,3,4).
\]
\end{proposition}
\begin{proof}
We identify the bosonizations and the corresponding super-data simultaneously.
For all four super-data representatives considered here,
$\gamma(x^ry^s)=(-1)^{r+s}$ and
$\operatorname{Ker}(\gamma)=\langle xy\rangle$.
We identify the abstract generator of $\Gamma=C_4$ with $xy$ in this kernel.
With this identification, the formulas of
Proposition~\ref{prp:induced-compatible-datum} give the four data above
without a further cochain change.
By Lemma~\ref{prp:Aut(H_c:sigma_0)}, the Hopf algebra
$H_{c:\sigma_0}$ has a unique super-datum up to equivalence, represented by
$(\kappa^{-1}(x^2y),\alpha_{y,-1})$.
The procedure of Section~\ref{sec:ell=0} associates $\DD_8^{(1)}$ with this representative.

By Lemma~\ref{prp:Aut(H_c:sigma_1)}, the three equivalence classes of super-data for $H_{c:\sigma_1}$ are represented by
\[
(\kappa^{-1}(x^2y),\alpha_{y,-1}),\quad
(\kappa^{-1}(x^2y),\alpha_{x^2y,1}),\quad
(\kappa^{-1}(x^2y),\alpha_{x^2y,-1}).
\]
Their associated compatible data are $\DD_8^{(2)}$, $\DD_8^{(3)}$, and $\DD_8^{(4)}$, respectively.
The four super-data classes are inequivalent, and the first also has a different bosonization from the other three.
Proposition~\ref{prp:one-to-one} therefore proves pairwise non-isomorphism.
\end{proof}

\paragraph{$\Gamma=C_2^2$ case}
Let $\Gamma=C_2^2=\langle y,z \mid y^2=z^2=1,zy=yz \rangle$.
Set
\[
\ttt(y) = z,\; \ttt(z) = y \quad\text{and}\quad
\ppp(1) = \ppp(yz) = 0, \; \ppp(y) = \ppp(z) = 1.
\]
We use the tuple notation $(\lambdaa)_{C_2^2}$ of Notation~\ref{not:cochain-values}.
It is known that $H^2(\Gamma) = H^2(C_2^2) \cong \mathbb{Z}_2$.
Writing elements as $y^rz^s$ with $r,s\in\{0,1\}$, we fix the normalized $2$-cocycle
\begin{equation}\label{eq:mu0}
\muu_0(y^rz^s,y^uz^v):=(-1)^{r(u+v)}
\qquad(r,s,u,v\in\{0,1\}).
\end{equation}
Since $\muu_0(y,z)=-1$ and $\muu_0(z,y)=1$, whereas every coboundary on an abelian group is symmetric, $\muu_0$ represents the non-trivial class in $H^2(\Gamma)$.

A direct computation shows that the following are compatible data for $\Gamma=C_2^2$.
\[
\begin{gathered}
\DD_8^{(5)}:=( \ttt, \ppp; \bigl(1,1,1,1\bigr), \muu_0 ), \quad
\DD_8^{(6)}:=( \ttt, \ppp; \bigl(1,\ii,\ii,1\bigr), 1 ), \\
\DD_8^{(7)}:=( \id, \ppp; \bigl(1,\ii,\ii,1\bigr), 1 ), \quad
\DD_8^{(8)}:=( \id, \ppp; \bigl(1,-\ii,-\ii,1\bigr), 1 ), \\
\DD_8^{(9)}:=( \id, \ppp; \bigl(1,\ii,\ii,1\bigr), \muu_0 ), \quad
\DD_8^{(10)}:=( \id, \ppp; \bigl(1,-\ii,\ii,-1\bigr), 1 ), \\
\DD_8^{(11)}:=( \id, \ppp; \bigl(1,\ii,-\ii,-1\bigr), \muu_0 ), \quad
\DD_8^{(12)}:=( \id, \ppp; \bigl(1,-\ii,\ii,-1\bigr), \muu_0 ).
\end{gathered}
\]
Set $\LL_8^{(i)}:=\LL(\Gamma,\DD_8^{(i)})$ for $i\in\{5,6,7,\dots,12\}$.
Note that $\G(\LL_8^{(i)}) \cong \G(\kk^\Gamma) \cong C_2^2$ since $\kk^\Gamma\cong \kk \Gamma$.

\begin{proposition} \label{prp:bosonization2}
For $i\in\{5,6,7,\dots,12\}$,
the bosonization of $\LL_8^{(i)}$ is isomorphic to $H_{d:1,1}$ if $i=5$, $H_{d:-1,1}$ if $i\in\{6,7,8,9\}$, and $H_{d:-1,-1}$ otherwise.
Moreover, the Hopf superalgebras $\LL_8^{(5)},\dots, \LL_8^{(12)}$ are pairwise non-isomorphic.
\end{proposition}
\begin{proof}
First, recall that $H_{d:1,1}$ has a unique super-datum $(\kappa^{-1}(zw),\alpha_{w,-1})$ (up to equivalence) by Lemma~\ref{prp:Aut(H_d:1,1)}.
The procedure in Section~\ref{sec:ell=0} associates the compatible datum
$\DD_8^{(5)}$ with this super-datum, up to equivalence.
Next, Lemma~\ref{prp:Aut(H_d:-1,1)} gives four super-data for $H_{d:-1,1}$ up to equivalence.
Under the construction of Proposition~\ref{prp:induced-compatible-datum}, their associated compatible data are, up to equivalence,
\[
\begin{aligned}
(\kappa^{-1}(zw),\alpha_{w,-1})&\longmapsto\DD_8^{(6)},\\
(\kappa^{-1}(zw),\alpha_{yz,-\ii})&\longmapsto\DD_8^{(8)},\\
(\kappa^{-1}(zw),\alpha_{yz,\ii})&\longmapsto\DD_8^{(7)},\\
(\kappa^{-1}(z),\alpha_{yz,-\ii})&\longmapsto\DD_8^{(9)}.
\end{aligned}
\]
Finally, the three super-data for $H_{d:-1,-1}$ correspond, up to equivalence, as follows:
\[
\begin{aligned}
(\kappa^{-1}(zw),\alpha_{yz,-\ii})&\longmapsto\DD_8^{(10)},\\
(\kappa^{-1}(z),\alpha_{yz,-\ii})&\longmapsto\DD_8^{(11)},\\
(\kappa^{-1}(z),\alpha_{yz,\ii})&\longmapsto\DD_8^{(12)}.
\end{aligned}
\]
Here are explicit choices for these equivalences.
For $\kappa(g)=zw$, identify the abstract generators $y,z$ of $\Gamma$
with the elements $y,zw$ in $G$;
for $\kappa(g)=z$, identify them with $y,yw$.
After these identifications, use the identity group automorphism and the
normalized cochains whose values at $(1,y,z,yz)$ are
\[
\begin{array}{c|c}
\text{compatible datum}&(f(1),f(y),f(z),f(yz))\\ \hline
\DD_8^{(5)}&(1,1,-1,-1)\\
\DD_8^{(6)},\DD_8^{(7)},\DD_8^{(8)},\DD_8^{(10)}&(1,\ii,1,-\ii)\\
\DD_8^{(9)},\DD_8^{(11)},\DD_8^{(12)}&(1,1,1,1).
\end{array}
\]
Substitution in the definition of equivalent compatible data gives exactly
the listed $\DD_8^{(5)},\dots,\DD_8^{(12)}$.
The super-data representatives are inequivalent, and the three bosonization types are distinguished as above; hence Proposition~\ref{prp:one-to-one} proves pairwise non-isomorphism.
\end{proof}

We have proved the following classification.
\begin{theorem}\label{thm:dim8-classification}
Let $\HH$ be a semisimple non-pointed Hopf superalgebra of dimension $8$ with $\HH_{\bar1}\neq0$.
Then $\HH$ is isomorphic to exactly one of $\LL_8^{(1)},\dots,\LL_8^{(12)}$.
\end{theorem}

Table~\ref{tab:dim8-classification} gives the groups of group-like elements and the bosonizations of these classes.
\begin{table}[H]
\centering
\caption{Bosonizations of the twelve semisimple classes of dimension $8$ with nonzero odd part}
\label{tab:dim8-classification}
\begin{tabular}{@{}ccc@{}}
\toprule
Classes & $\G(\HH)$ & $\widehat{\HH}$\\
\midrule
$\LL_8^{(1)}$ & $C_4$ & $H_{c:\sigma_0}$\\
$\LL_8^{(2)},\LL_8^{(3)},\LL_8^{(4)}$ & $C_4$ & $H_{c:\sigma_1}\cong(H_{d:-1,-1})^*$\\
$\LL_8^{(5)}$ & $C_2^2$ & $H_{d:1,1}$\\
$\LL_8^{(6)},\LL_8^{(7)},\LL_8^{(8)},\LL_8^{(9)}$ & $C_2^2$ & $H_{d:-1,1}$\\
$\LL_8^{(10)},\LL_8^{(11)},\LL_8^{(12)}$ & $C_2^2$ & $H_{d:-1,-1}$\\
\bottomrule
\end{tabular}
\end{table}

\begin{remark*}
Table~\ref{tab:dim8-classification} shows that pairwise non-isomorphic Hopf superalgebras can have the same bosonization.
For example, $\LL_8^{(2)}$, $\LL_8^{(3)}$, and $\LL_8^{(4)}$ all have bosonization $H_{c:\sigma_1}$.
By the tensor equivalence recalled in Subsection~\ref{subsec:Hopf-super}, their categories of finite-dimensional supermodules with even morphisms are equivalent as ordinary tensor categories.
Thus, the classification distinguishes super-forms that are not detected by the bosonized Hopf algebra alone.
\end{remark*}

\subsubsection{Identification with known objects}
The construction of Theorem~\ref{thm:AEG} can also be applied directly to two familiar semisimple Hopf algebras.
With our dihedral convention, let
\[
D_8=\langle a,b\mid a^2=b^4=1,\ aba=b^{-1}\rangle.
\]
The relevant pair is $(\kk D_8,a)$.
The second pair is $(\AKP,s)$.
For the Kac--Paljutkin Hopf algebra $\AKP$, we use the following presentation, which is equivalent to the abelian-extension realization in Section~\ref{sec:dim=4}; see~\cite{Mas95b}:
\[
\begin{gathered}
\AKP=
\bigangle{r,s,q}{r^2=s^2=1,\ q^2=\frac12(1+r+s-rs),\\
rs=sr,\ qr=sq,\ rq=qs},\\
\Delta(r)=r\otimes r,\qquad \Delta(s)=s\otimes s,\\
\Delta(q)=\frac12\bigl(q\otimes q+q\otimes rq+sq\otimes q-sq\otimes rq\bigr).
\end{gathered}
\]
Here $r=\kappa^{-1}(x)$, $s=\kappa^{-1}(y)$, and $q=t$ under the abelian-extension realization of Section~\ref{sec:dim=4}.
We now identify the corresponding members of our list.

\begin{proposition}\label{prp:K5-D8}
We have $\LL_8^{(5)}\cong(\kk D_8)^a$.
\end{proposition}
\begin{proof}
Let $\Gamma=C_2^2=\langle y,z\mid y^2=z^2=1,\ yz=zy\rangle$, and let $\vartheta\in\Aut(\Gamma)$ interchange $y$ and $z$.
Put
\[
E:=\langle b^2,a\rangle\cong C_2^2.
\]
For $g=y^rz^s\in\Gamma$, define $\chi_g\in\widehat E$ by
\[
\chi_g(b^2)=(-1)^{r+s},
\qquad
\chi_g(a)=(-1)^s,
\]
and put
\[
P_g:=\frac14\sum_{e\in E}\chi_g(e)e,
\qquad
Q_g:=(-1)^rP_gb.
\]
We regard these elements as elements of $(\kk D_8)^a$, whose underlying algebra is $\kk D_8$.
The elements $P_g$ are pairwise orthogonal primitive idempotents of $\kk E$, and their sum is $1$.
The relations in $D_8$ give
\[
bP_gb^{-1}=P_{\vartheta(g)},
\qquad
P_gb^2=(-1)^{r+s}P_g,
\qquad
aQ_ga=(-1)^{r+s}Q_g.
\]
Consequently, with $\ppp(y^rz^s)=r+s\pmod2$, the $P_g$ are even, $|Q_g|=\ppp(g)$, and
\begin{equation}\label{eq:K5-mult}
\begin{aligned}
P_gP_h&=\delta_{g,h}P_g,&
P_gQ_h&=\delta_{g,h}Q_g,\\
Q_gP_h&=\delta_{\vartheta(g),h}Q_g,&
Q_gQ_h&=\delta_{\vartheta(g),h}P_g.
\end{aligned}
\end{equation}

The coproduct $\Delta$ of $\kk D_8$ is
\[
\Delta(P_g)=\sum_{uv=g}P_u\otimes P_v,
\qquad
\Delta(Q_g)=\sum_{uv=g}Q_u\otimes Q_v.
\]
Recall from Theorem~\ref{thm:AEG} that, for a homogeneous element $h$,
\[
\Delta^a(h)=\Delta_{\bar0}(h)-(-1)^{|h|}(a\otimes1)\Delta_{\bar1}(h),
\]
where $\Delta_{\bar\epsilon}(h)$ is the part whose second tensor factor has parity $\bar\epsilon$.
Write $u=y^rz^s$.
In the summand $Q_u\otimes Q_v$, the second factor has parity $\ppp(v)$.
If $\ppp(v)=0$, this summand is unchanged.
If $\ppp(v)=1$, the formula above and $aQ_u=(-1)^sQ_u$ multiply it by $-(-1)^{\ppp(g)+s}$.
Since $g=uv$, in both cases the coefficient is
\[
(-1)^{r\ppp(v)}=\muu_0(u,v)
\]
by~\eqref{eq:mu0}. Hence
\begin{equation}\label{eq:K5-comult}
\Delta^a(P_g)=\sum_{uv=g}P_u\otimes P_v,
\qquad
\Delta^a(Q_g)=\sum_{uv=g}\muu_0(u,v)Q_u\otimes Q_v.
\end{equation}
Equations~\eqref{eq:K5-mult} and~\eqref{eq:K5-comult}, together with the counit, are precisely the defining formulas for $\LL(\Gamma,\DD_8^{(5)})$.
Thus
\[
v(g,\bar0)\longmapsto P_g,
\qquad
v(g,\bar1)\longmapsto Q_g
\]
defines the desired Hopf superalgebra isomorphism.
\end{proof}

\begin{proposition}\label{prp:K6-KP}
We have $\LL_8^{(6)}\cong(\AKP)^s$.
\end{proposition}
\begin{proof}
Again write $\Gamma=C_2^2=\langle y,z\mid y^2=z^2=1,\ yz=zy\rangle$ and let $\vartheta$ interchange $y$ and $z$.
Put
\[
E':=\langle r,s\rangle\cong C_2^2.
\]
For $g=y^mz^n$, define $\chi_g\in\widehat{E'}$ by
\[
\chi_g(r)=(-1)^n,
\qquad
\chi_g(s)=(-1)^m,
\qquad
P_g:=\frac14\sum_{e\in E'}\chi_g(e)e.
\]
The elements $P_g$ are pairwise orthogonal primitive idempotents of $\kk E'$, and their sum is $1$.
The defining relations of $\AKP$ imply
\[
qP_g=P_{\vartheta(g)}q,
\qquad
P_gq^2=(-1)^{mn}P_g.
\]
Put
\[
\rho(g):=\ii^n(-1)^{mn},
\qquad
Q_g:=\rho(g)P_gq,
\qquad
\lambdaa_6(g):=\ii^{m+n}(-1)^{mn}.
\]
We will show that $\Delta^s(Q_g)=\sum_{uv=g}Q_u\otimes Q_v$.
We regard $P_g$ and $Q_g$ as elements of $(\AKP)^s$, whose underlying algebra is $\AKP$.
Conjugation by $s$ shows that $P_g$ is even and $|Q_g|=\ppp(g)=m+n\pmod2$.
Moreover,
\begin{equation}\label{eq:K6-mult}
\begin{aligned}
P_gP_h&=\delta_{g,h}P_g,&
P_gQ_h&=\delta_{g,h}Q_g,\\
Q_gP_h&=\delta_{\vartheta(g),h}Q_g,&
Q_gQ_h&=\delta_{\vartheta(g),h}\lambdaa_6(g)P_g.
\end{aligned}
\end{equation}
Using Notation~\ref{not:cochain-values}, we have
\[
(\lambdaa_6)_{C_2^2}=\bigl(1,\ii,\ii,1\bigr).
\]

Expanding the formula for $\Delta(q)$ gives
\[
\Delta(P_g)=\sum_{uv=g}P_u\otimes P_v,
\qquad
\Delta(P_gq)=\sum_{uv=g}(-1)^{r_us_v}P_uq\otimes P_vq,
\]
where $u=y^{r_u}z^{s_u}$ and $v=y^{r_v}z^{s_v}$.
The definition of $\rho$ gives
\begin{equation}\label{eq:rho-identity}
\rho(g)\rho(u)^{-1}\rho(v)^{-1}(-1)^{r_us_v}
 =(-1)^{s_u\ppp(v)}.
\end{equation}
Since $sQ_u=(-1)^{r_u}Q_u$, the formula in Theorem~\ref{thm:AEG} multiplies the coefficient of $Q_u\otimes Q_v$ by $(-1)^{s_u\ppp(v)}$.
Together with~\eqref{eq:rho-identity}, this gives a coefficient of $1$. Hence
\begin{equation}\label{eq:K6-comult}
\Delta^s(P_g)=\sum_{uv=g}P_u\otimes P_v,
\qquad
\Delta^s(Q_g)=\sum_{uv=g}Q_u\otimes Q_v.
\end{equation}
Together with the counit, equations~\eqref{eq:K6-mult} and~\eqref{eq:K6-comult} are exactly the structure formulas for $\LL(\Gamma,\DD_8^{(6)})$.
Thus $v(g,\bar0)\mapsto P_g$ and $v(g,\bar1)\mapsto Q_g$ give an isomorphism $\LL_8^{(6)}\cong(\AKP)^s$.
\end{proof}

The tensor product $\LL_4^{(i)}\otimes\kk C_2$ ($i\in\{1,2\}$) is semisimple of dimension $8$ and has nonzero odd part.
It is therefore non-pointed, since a semisimple pointed Hopf superalgebra is purely even.
Theorem~\ref{thm:dim8-classification} implies that it is isomorphic to one of $\LL_8^{(j)}$ ($j\in\{1,\dots,12\}$).
The next proposition determines which one it is.
\begin{proposition} \label{prp:Masuoka-split}
We have $\LL_8^{(7)}\cong \LL_4^{(1)}\otimes \kk C_2$ and $\LL_8^{(8)}\cong \LL_4^{(2)}\otimes \kk C_2$.
\end{proposition}
\begin{proof}
For $\delta\in\{1,-1\}$, let
\[
C_2=\langle x\mid x^2=1\rangle,
\qquad
C_2^2=\langle y,z\mid y^2=z^2=1,\ yz=zy\rangle.
\]
Let $\ppp$ denote the non-trivial homomorphism on $C_2$ and the homomorphism on $C_2^2$ determined by
\[
\ppp(x)=1,
\qquad
\ppp(y)=\ppp(z)=1.
\]
Define cochains $\lambdaa_\delta:C_2\to\kk^\times$ and $\lambdaa'_\delta:C_2^2\to\kk^\times$ by
\[
\lambdaa_\delta(x)=\delta\ii,
\qquad
\lambdaa'_\delta(y)=\lambdaa'_\delta(z)=\delta\ii,
\qquad
\lambdaa'_\delta(yz)=1,
\]
and put
\[
\mathscr U_\delta:=\LL(C_2,(\id,\ppp;\lambdaa_\delta,1)),
\qquad
\mathscr V_\delta:=\LL(C_2^2,(\id,\ppp;\lambdaa'_\delta,1)).
\]
Then
\[
\mathscr U_1=\LL_4^{(1)},
\qquad
\mathscr U_{-1}=\LL_4^{(2)},
\qquad
\mathscr V_1=\LL_8^{(7)},
\qquad
\mathscr V_{-1}=\LL_8^{(8)}.
\]
By the definition of $\LL(\Gamma,\DD)$, write
$u_\delta(g,\bar\epsilon)$ for the basis vector of $\mathscr U_\delta$ indexed by $(g,\bar\epsilon)\in C_2\times\Z_2$, and
$v_\delta(h,\bar\epsilon)$ for the basis vector of $\mathscr V_\delta$ indexed by $(h,\bar\epsilon)\in C_2^2\times\Z_2$.
Let $C_2=\langle t\rangle$ in the second tensor factor.
For $\epsilon\in\{0,1\}$, define
\[
\begin{aligned}
v_\delta(1,\bar\epsilon)&\longmapsto
 \frac12 u_\delta(1,\bar\epsilon)\otimes(1+t),&
 v_\delta(y,\bar\epsilon)&\longmapsto
 \frac12 u_\delta(x,\bar\epsilon)\otimes(1+t),\\
v_\delta(z,\bar\epsilon)&\longmapsto
 \frac12 u_\delta(x,\bar\epsilon)\otimes(1-t),&
 v_\delta(yz,\bar\epsilon)&\longmapsto
 \frac12 u_\delta(1,\bar\epsilon)\otimes(1-t).
\end{aligned}
\]
The elements $(1+t)/2$ and $(1-t)/2$ are orthogonal idempotents whose sum is $1$.
Substitution in the defining formulas of $\LL(\Gamma,\DD)$ shows that this assignment preserves multiplication, comultiplication, unit, and counit.
Its images form a basis, so it is an isomorphism
$\mathscr V_\delta\cong\mathscr U_\delta\otimes\kk C_2$.
Taking $\delta=1$ and $\delta=-1$ proves the two assertions.
\end{proof}

\section{Duals of semisimple non-pointed Hopf superalgebras of dimension \texorpdfstring{$8$}{8}}\label{sec:dual}
In this section, we determine the duals of the semisimple non-pointed Hopf superalgebras of dimension $8$ classified in Section~\ref{sec:class}.
As in Subsection~\ref{subsec:Hopf-super}, we use the tensor-pairing convention
$\langle f\otimes g,a\otimes b\rangle=\langle f,a\rangle\langle g,b\rangle$.

For a basis $\{v_1,\dots,v_n\}$ of a finite-dimensional vector space $V$,
we let $\{v_1^*,\dots,v_n^*\}$ denote its dual basis of $V^*$.
We first prove the following lemma.
\begin{lemma} \label{prp:dual-lemma}
Take $\Gamma,\Gamma'\in\{C_2^2,C_4\}$ and suppose that $\LL(\Gamma,\DD)^*\cong\LL(\Gamma',\DD')$ as Hopf superalgebras
for some $\DD=(\ttt,\ppp;\lambdaa,\muu)\in\CD(\Gamma)$ and $\DD'=(\ttt',\ppp';\lambdaa',\muu') \in \CD(\Gamma')$.
If $\muu=1$, then $\ttt'=\id$.
\end{lemma}
\begin{proof}
Since $\muu=1$, for all $g,h\in\Gamma$ and $\epsilon,\eta\in\{0,1\}$ we have
\[
v(g,\bar\epsilon)^*v(h,\bar\eta)^*
 =\delta_{\epsilon,\eta}v(gh,\bar\epsilon)^*
.
\]
Thus the underlying ordinary algebra of $\LL(\Gamma,\DD)^*$ is isomorphic to $\kk\Gamma\times\kk\Gamma$ and is commutative, since $\Gamma$ is abelian.
If $\ttt'(g)\neq g$ for some $g\in\Gamma'$, then in $\LL(\Gamma',\DD')$ we have
\[
v(g,\bar0)v(g,\bar1)=v(g,\bar1),
\qquad
v(g,\bar1)v(g,\bar0)=0,
\]
contradicting commutativity.
Therefore $\ttt'=\id$.
\end{proof}

The same argument applies when $\Gamma$ is any finite abelian group and $\Gamma'$ is any finite group.

We next give two explicit duality isomorphisms.
\begin{lemma}\label{lem:explicit-duality}
For $\delta\in\{1,-1\}$, let
\[
\mathscr C_\delta:=\LL(C_4,(\iota,\ppp;\lambdaa_\delta,1)),
\qquad
\mathscr E_\delta:=\LL(C_2^2,(\id,\ppp;\lambdaa'_\delta,\muu_0)),
\]
where
\[
C_4=\langle x\mid x^4=1\rangle,
\qquad
C_2^2=\langle y,z\mid y^2=z^2=1,\ yz=zy\rangle,
\]
$\iota(x)=x^{-1}$, $\ppp(x)=1$, and $\ppp(y)=\ppp(z)=1$.
Define $\lambdaa_\delta$ by
\[
\lambdaa_\delta(x)=\lambdaa_\delta(x^3)=\delta\ii,
\qquad
\lambdaa_\delta(x^2)=1,
\]
and, using Notation~\ref{not:cochain-values}, define $\lambdaa'_\delta$ by
\[
(\lambdaa'_\delta)_{C_2^2}=
\bigl(1,-\delta\ii,\delta\ii,-1\bigr).
\]
Then $\mathscr E_\delta\cong\mathscr C_\delta^*$.
In particular,
\[
\LL_8^{(12)}\cong(\LL_8^{(3)})^*,
\qquad
\LL_8^{(11)}\cong(\LL_8^{(4)})^*.
\]
\end{lemma}
\begin{proof}
Choose $\zeta_\delta\in\kk$ with $\zeta_\delta^2=-\delta\ii$.
By the definition of $\LL(\Gamma,\DD)$, let
$u_\delta(g,\bar\epsilon)$ be the basis vector of $\mathscr E_\delta$ indexed by $(g,\bar\epsilon)\in C_2^2\times\Z_2$, and let
$w_\delta(x^j,\bar\epsilon)$ be the basis vector of $\mathscr C_\delta$ indexed by $(x^j,\bar\epsilon)\in C_4\times\Z_2$.
Put $e_{j,\bar\epsilon}:=w_\delta(x^j,\bar\epsilon)^*$, with $j$ read modulo $4$.
Define $\Phi_\delta:\mathscr E_\delta\to\mathscr C_\delta^*$ by Table~\ref{tab:explicit-duality}.

\begin{table}[H]
\centering
\caption{The explicit duality isomorphism $\Phi_\delta:\mathscr E_\delta\to\mathscr C_\delta^*$}
\label{tab:explicit-duality}
\begin{tabular}{@{}c c c@{}}
\toprule
$g$ & $\Phi_\delta(u_\delta(g,\bar0))$ & $\Phi_\delta(u_\delta(g,\bar1))$\\
\addlinespace[3pt]
\midrule
$1$  & $\frac12(e_{0,\bar0}+e_{2,\bar0})$
   & $\frac12(e_{1,\bar0}+e_{3,\bar0})$\\
\addlinespace[3pt]
$y$  & $\frac12(e_{0,\bar1}+e_{2,\bar1})$
   & $\frac{\zeta_\delta}{2}(e_{1,\bar1}+e_{3,\bar1})$\\
\addlinespace[3pt]
$z$  & $\frac12(e_{0,\bar1}-e_{2,\bar1})$
   & $\frac{\zeta_\delta}{2}(e_{1,\bar1}-e_{3,\bar1})$\\
\addlinespace[3pt]
$yz$ & $\frac12(e_{0,\bar0}-e_{2,\bar0})$
   & $\frac12(-e_{1,\bar0}+e_{3,\bar0})$\\
\addlinespace[3pt]
\bottomrule
\end{tabular}
\end{table}

Since the cocycle of $\mathscr C_\delta$ is trivial, dualizing its defining formulas gives
\begin{equation}\label{eq:dual-Cdelta}
\begin{aligned}
e_{j,\bar\epsilon}e_{k,\bar\eta}
 &=\delta_{\epsilon,\eta}\,e_{j+k,\bar\epsilon},\\
\Delta(e_{j,\bar\rho})
 &=\sum_{\epsilon+\eta=\rho}\lambdaa_\delta(x^j)^{\epsilon\eta}
 e_{j,\bar\epsilon}\otimes e_{(-1)^\epsilon j,\bar\eta},
\end{aligned}
\end{equation}
where all subscripts are taken modulo $4$ or $2$, as appropriate.
The eight vectors in Table~\ref{tab:explicit-duality} form a homogeneous basis of $\mathscr C_\delta^*$.
For example, the row indexed by $y$ gives
\[
\left(\frac{\zeta_\delta}{2}(e_{1,\bar1}+e_{3,\bar1})\right)^2
 =\frac{\zeta_\delta^2}{2}(e_{0,\bar1}+e_{2,\bar1})
 =-\delta\ii\,\Phi_\delta(u_\delta(y,\bar0)),
\]
which explains the condition $\zeta_\delta^2=-\delta\ii$.
Substitution into~\eqref{eq:dual-Cdelta}, using
$\zeta_\delta^2=-\delta\ii$ and
$\muu_0(y^rz^\ell,y^uz^v)=(-1)^{r(u+v)}$, gives exactly the multiplication and comultiplication of $\mathscr E_\delta$; the unit and counit are immediate.
Thus $\Phi_\delta$ is a Hopf superalgebra isomorphism.
For $\delta=1$, $(\mathscr C_\delta,\mathscr E_\delta)=(\LL_8^{(3)},\LL_8^{(12)})$, while for $\delta=-1$, $(\mathscr C_\delta,\mathscr E_\delta)=(\LL_8^{(4)},\LL_8^{(11)})$.
\end{proof}

Proposition~\ref{prp:Masuoka-split} and the self-duality of $\LL_4^{(1)}$ and $\LL_4^{(2)}$ imply that $\LL_8^{(7)}$ and $\LL_8^{(8)}$ are self-dual.
Propositions~\ref{prp:dual-hat},~\ref{prp:bosonization1}, and~\ref{prp:bosonization2} determine some of these duals from the duals of their bosonizations.
By~\cite[Remark~1.4]{Kas00}, the Hopf algebras $H_{c:\sigma_0}$, $H_{d:1,1}$, and $H_{d:-1,1}$ are self-dual.
The direct isomorphism in Lemma~\ref{lem:explicit-duality}, together with Propositions~\ref{prp:dual-hat},~\ref{prp:bosonization1}, and~\ref{prp:bosonization2}, gives
\[
H_{c:\sigma_1}^*
\cong \widehat{(\LL_8^{(3)})^*}
\cong \widehat{\LL_8^{(12)}}
\cong H_{d:-1,-1}.
\]
Kashina's classification identifies the two duals $(H_{c:\sigma_1})^*$ and $(H_{b:1})^*$, in some order, with $H_{d:-1,-1}$ and $H_{d:1,-1}$~\cite[Section~9]{Kas00}; hence $(H_{b:1})^*\cong H_{d:1,-1}$ as well.
The bosonizations $H_{c:\sigma_0}$ and $H_{d:1,1}$ each have a unique non-trivial super-form up to isomorphism; Proposition~\ref{prp:dual-hat} therefore implies that $\LL_8^{(1)}$ and $\LL_8^{(5)}$ are self-dual.
The dual of $\LL_8^{(6)}$ must be one of $\LL_8^{(6)},\LL_8^{(7)},\LL_8^{(8)},\LL_8^{(9)}$, since these are precisely the non-trivial super-forms of $H_{d:-1,1}$.
Duality is an involution on isomorphism classes, so the already self-dual classes $\LL_8^{(7)}$ and $\LL_8^{(8)}$ cannot occur.
Lemma~\ref{prp:dual-lemma} excludes $\LL_8^{(6)}$ itself.
Hence $(\LL_8^{(6)})^*\cong\LL_8^{(9)}$.
Consequently, up to isomorphism, the three duals on the left are precisely the three classes on the right, in some order:
\[
\{(\LL_8^{(2)})^*,(\LL_8^{(3)})^*,(\LL_8^{(4)})^*\}
=
\{\LL_8^{(10)},\LL_8^{(11)},\LL_8^{(12)}\}.
\]
Lemma~\ref{lem:explicit-duality} identifies two of these pairs, and the third follows by elimination.

We have proved the following theorem.
\begin{theorem}
For $i\in\{1,5,7,8\}$, $\LL_8^{(i)}$ is self-dual.
The duals of $\LL_8^{(2)}$, $\LL_8^{(3)}$, $\LL_8^{(4)}$, and $\LL_8^{(6)}$ are $\LL_8^{(10)}$, $\LL_8^{(12)}$, $\LL_8^{(11)}$, and $\LL_8^{(9)}$, respectively.
\end{theorem}

\appendix
\section{Complete classification in dimensions up to \texorpdfstring{$10$}{10}} \label{sec:appendix}

This appendix summarizes the classification of Hopf superalgebras of dimension at most $10$ obtained in~\cite{SSW24,SW25} and in the present paper.
We first give the number of isomorphism classes in each dimension and then list every class with nonzero odd part.
Purely even classes are included in the counts, but their presentations are not repeated.

\subsection{Numbers of isomorphism classes}\label{subsec:app-overview}
For a non-semisimple Hopf superalgebra $\mathcal H$ with $\mathcal H_{\bar1}\neq0$, we distinguish the following three mutually exclusive cases:
\begin{enumerate}
\item $\mathcal H$ is pointed;
\item $\mathcal H$ is not pointed but $\mathcal H^*$ is pointed;
\item neither $\mathcal H$ nor $\mathcal H^*$ is pointed.
\end{enumerate}
Together with the purely even and the non-trivial semisimple cases, these cases exhaust all Hopf superalgebras of dimension at most $10$.

By Theorem~\ref{thm:main1}, every Hopf superalgebra of dimension $9$ is purely even.
Table~\ref{tab:app-counts} gives the resulting numbers.

\begin{table}[H]
\centering
\caption{Numbers of isomorphism classes of Hopf superalgebras of dimension at most $10$}
\label{tab:app-counts}
\begin{tabular}{@{}c r r r r r r@{}}
\toprule
$\dim(\mathcal H)$ & Purely even & \begin{tabular}[c]{@{}c@{}}Semisimple\\$\mathcal H_{\bar1}\neq0$\end{tabular}
& \multicolumn{3}{c}{Non-semisimple with $\mathcal H_{\bar1}\neq0$} & Total\\
\cmidrule(lr){4-6}
& & & Pointed & \begin{tabular}[c]{@{}c@{}}Not pointed;\\pointed dual\end{tabular} & Neither & \\
\midrule
1  & 1  & 0  & 0  & 0 & 0 & 1  \\
2  & 1  & 0  & 1  & 0 & 0 & 2  \\
3  & 1  & 0  & 0  & 0 & 0 & 1  \\
4  & 3  & 2  & 4  & 0 & 0 & 9  \\
5  & 1  & 0  & 0  & 0 & 0 & 1  \\
6  & 3  & 1  & 4  & 1 & 0 & 9  \\
7  & 1  & 0  & 0  & 0 & 0 & 1  \\
8  & 14 & 12 & 19 & 2 & 8 & 55 \\
9  & 4  & 0  & 0  & 0 & 0 & 4  \\
10 & 3  & 1  & 4  & 1 & 0 & 9  \\
\midrule
Total & 32 & 16 & 32 & 4 & 8 & 92 \\
\bottomrule
\end{tabular}
\end{table}

The numbers in the purely even column follow from the classification of ordinary Hopf algebras.
In dimensions $1,2,3,5$, and $7$ there is one class.
In dimension $4$ the three classes are $\kk C_4$, $\kk C_2^2$, and the Sweedler
Hopf algebra.
In dimension $6$ they are $\kk C_6$, $\kk\mathfrak{S}_3$, and $\kk^{\mathfrak{S}_3}$, and in
dimension $10$ they are $\kk C_{10}$, $\kk D_{10}$, and $\kk^{D_{10}}$.
In dimension $8$ there are fourteen ordinary classes: eight semisimple and six
non-semisimple classes.
In dimension $9$ the four classes are $\kk C_9$, $\kk C_3^2$, and the two Taft
algebras corresponding to the two primitive cube roots of unity; see~\cite[Theorem~6.5]{Ng02}.
See, for example, \cite{BG13,Ste99} and the references therein.

\subsection{Notation}
Table~\ref{tab:app-notation} lists the notation used in this appendix and the references for the corresponding classifications.

\begin{table}[H]
\centering
\caption{Notation for the non-trivial Hopf superalgebras}
\label{tab:app-notation}
\begin{tabularx}{\textwidth}{@{}>{\raggedright\arraybackslash}X c >{\raggedright\arraybackslash}X l@{}}
\toprule
Notation & Dimensions & Type & References \\
\midrule
$\Lambda_1$, $\HH_n^{(i)}$ & $n=2,4,6,8,10$ & non-semisimple pointed & \cite{SSW24} \\
$(\HH_6^{(4)})^*$, $(\HH_8^{(18)})^*$, $(\HH_8^{(19)})^*$, $(\HH_{10}^{(4)})^*$
& $6,8,10$ & non-semisimple; not pointed, with pointed dual & \cite{SSW24} \\
$K_8(\zeta;\epsilon,\eta)$ & $8$ & non-semisimple; not pointed, with non-pointed dual & \cite{SW25} \\
$\LL_4^{(i)}$, $\LL_8^{(i)}$, $(\kk D_{2p})^a$
& $4,6,8,10$ & semisimple with nonzero odd part & \cite{SW25} and this paper \\
\bottomrule
\end{tabularx}
\end{table}

\subsection{Non-semisimple pointed cases}
We first recall the notation of \cite{SSW24}.
Let $\Gamma$ be a finite abelian group and let $\theta$ be a positive integer.
Let
\[
 \DD_{\mathrm{pt}}
 =\bigl((g_i,\chi_i,\nu_i;\epsilon_i)\bigr)_{i=1}^{\theta}
 \in\bigl(\Gamma\times\widehat\Gamma\times\{0,1\}\times\{0,1\}\bigr)^\theta
\]
be a compatible datum in the sense of that paper.
Here the symbol $\nu_i$ is used for the third entry to avoid confusion with the
$2$-cocycle $\muu$ used in the present paper; it is denoted by
$\mu_i$ in \cite{SSW24}.
Put
\[
 N_i=\ord\bigl((-1)^{\epsilon_i}\chi_i(g_i)\bigr)>1.
\]
The pointed Hopf superalgebra $\AA(\Gamma,\DD_{\mathrm{pt}})$ is generated by
$\Gamma$ and homogeneous elements $z_1,\dots,z_\theta$ of degrees
$|z_i|=\epsilon_i$, with
\[
 hz_i=\chi_i(h)z_ih,\qquad
 z_i^{N_i}=\nu_i(1-g_i^{N_i}),\qquad
 z_iz_j=(-1)^{\epsilon_i\epsilon_j}\chi_j(g_i)z_jz_i\quad(i\neq j),
\]
and
\[
 \Delta(z_i)=g_i\otimes z_i+z_i\otimes1.
\]
We write $\Lambda_r:=\bigwedge(z_1,\dots,z_r)$ for the exterior Hopf
superalgebra on $r$ odd primitive generators; thus
\[
 \Lambda_r=\AA\bigl(\{1\},((1,1,0;1))_{i=1}^r\bigr).
\]

\subsubsection{Dimension \texorpdfstring{$2$}{2}}
There is exactly one non-trivial class, namely the exterior Hopf superalgebra $\Lambda_1$. It is self-dual.

\subsubsection{Dimensions \texorpdfstring{$2p$}{2p} with \texorpdfstring{$p=3,5$}{p=3,5}}
Fix a primitive $p$-th root of unity $\zeta_p\in\kk$.
Let $C_p=\langle g\mid g^p=1\rangle$, let $\chi:C_p\to\kk^\times$ be the character determined by $\chi(g)=\zeta_p$, and put $h=g^{(p+1)/2}$.
For each $p\in\{3,5\}$, the four classes are listed in Table~\ref{tab:app-pointed-2p}.

\begin{table}[H]
\centering
\caption{Non-semisimple pointed classes of dimension $2p$, $p\in\{3,5\}$}
\label{tab:app-pointed-2p}
\begin{tabularx}{\textwidth}{@{}c X c@{}}
\toprule
Class & Realization & Dual \\
\midrule
$\HH_{2p}^{(1)}$ & $\kk C_p\otimes\Lambda_1=\AA(C_p,(1,1,0;1))$ & self-dual \\
$\HH_{2p}^{(2)}$ & $\AA(C_p,(h,1,0;1))$ & $\HH_{2p}^{(3)}$ \\
$\HH_{2p}^{(3)}$ & $\AA(C_p,(1,\chi,0;1))$ & $\HH_{2p}^{(2)}$ \\
$\HH_{2p}^{(4)}$ & $\AA(C_p,(h,1,1;1))$ & non-pointed \\
\bottomrule
\end{tabularx}
\end{table}
The duality statements follow from~\cite[Theorem~5.1]{SSW24}.
The non-pointed dual of the last row is listed separately in
Subsection~\ref{subsec:app-dual-pointed}.

\subsubsection{Dimension \texorpdfstring{$4$}{4}}
Let $C_2=\langle g\mid g^2=1\rangle$ and let $\chi(g)=-1$. The four classes are listed in Table~\ref{tab:app-pointed-4}.

\begin{table}[H]
\centering
\caption{Non-semisimple pointed classes of dimension $4$}
\label{tab:app-pointed-4}
\begin{tabularx}{\textwidth}{@{}c X c@{}}
\toprule
Class & Realization & Dual \\
\midrule
$\HH_4^{(1)}$ & $\Lambda_2$ & self-dual \\
$\HH_4^{(2)}$ & $\kk C_2\otimes\Lambda_1=\AA(C_2,(1,1,0;1))$ & self-dual \\
$\HH_4^{(3)}$ & $\AA(C_2,(g,1,0;1))$ & $\HH_4^{(4)}$ \\
$\HH_4^{(4)}$ & $\AA(C_2,(1,\chi,0;1))$ & $\HH_4^{(3)}$ \\
\bottomrule
\end{tabularx}
\end{table}

\subsubsection{Dimension \texorpdfstring{$8$}{8}}
There are nineteen classes, whose duals are given in~\cite[Theorem~5.3]{SSW24}.
We group the classes by their group $\Gamma$ of group-like elements.
For $\Gamma=C_2=\langle g\rangle$, let $\chi(g)=-1$.
The notation $T_4(-1)$ denotes the four-dimensional Sweedler Hopf algebra.

The class $\HH_8^{(19)}$ is not of the form $\AA(\Gamma,\DD_{\mathrm{pt}})$ above.
By~\cite[Section~4.4.6]{SSW24}, it is generated by an even group-like element $g$ and odd elements $x,y$ with relations
\begin{equation}\label{eq:app-H19}
\begin{gathered}
g^2=1,\qquad gx=xg,\qquad gy=yg,\\
x^2=y^2=0,\qquad xy+yx=1-g,
\end{gathered}
\end{equation}
and coproduct
\[
\Delta(x)=1\otimes x+x\otimes1,
\qquad
\Delta(y)=g\otimes y+y\otimes1.
\]

The classes with $\Gamma=\{1\}$ or $C_2$ are listed in Table~\ref{tab:app-pointed-8-C2}.

\begin{table}[H]
\centering
\caption{Non-semisimple pointed classes of dimension $8$: $\Gamma=\{1\}$ or $C_2$}
\label{tab:app-pointed-8-C2}
\begin{tabularx}{\textwidth}{@{}c c X c@{}}
\toprule
Class & $\Gamma$ & Datum $\DD_{\mathrm{pt}}$ or realization & Dual \\
\midrule
$\HH_8^{(1)}$ & $\{1\}$ & $\Lambda_3$ & self-dual \\
$\HH_8^{(2)}$ & $C_2$ & $\kk C_2\otimes\Lambda_2$ & self-dual \\
$\HH_8^{(3)}$ & $C_2$ & $((1,1,0;1),(g,1,0;1))$ & $\HH_8^{(4)}$ \\
$\HH_8^{(4)}$ & $C_2$ & $((1,1,0;1),(1,\chi,0;1))$ & $\HH_8^{(3)}$ \\
$\HH_8^{(5)}$ & $C_2$ & $((g,1,0;1),(g,1,0;1))$ & $\HH_8^{(6)}$ \\
$\HH_8^{(6)}$ & $C_2$ & $((1,\chi,0;1),(1,\chi,0;1))$ & $\HH_8^{(5)}$ \\
$\HH_8^{(7)}$ & $C_2$ & $T_4(-1)\otimes\Lambda_1
 =\AA(C_2,((g,\chi,0;0),(1,1,0;1)))$ & self-dual \\
$\HH_8^{(19)}$ & $C_2$ & \eqref{eq:app-H19} & non-pointed \\
\bottomrule
\end{tabularx}
\end{table}

For $\Gamma=C_2^2=\langle g_1,g_2\rangle$, let
\[
 \chi_1(g_1)=\chi_2(g_2)=-1,
 \qquad
 \chi_1(g_2)=\chi_2(g_1)=1.
\]
The corresponding classes are listed in Table~\ref{tab:app-pointed-8-C22}.

\begin{table}[H]
\centering
\caption{Non-semisimple pointed classes of dimension $8$: $\Gamma=C_2^2$}
\label{tab:app-pointed-8-C22}
\begin{tabularx}{\textwidth}{@{}c X c@{}}
\toprule
Class & Datum $\DD_{\mathrm{pt}}$ or realization & Dual \\
\midrule
$\HH_8^{(8)}$ & $\kk C_2^2\otimes\Lambda_1=\AA(C_2^2,(1,1,0;1))$ & self-dual \\
$\HH_8^{(9)}$ & $\AA(C_2^2,(g_1,1,0;1))$ & $\HH_8^{(10)}$ \\
$\HH_8^{(10)}$ & $\AA(C_2^2,(1,\chi_1,0;1))$ & $\HH_8^{(9)}$ \\
$\HH_8^{(11)}$ & $\AA(C_2^2,(g_1,\chi_2,0;1))$ & self-dual \\
\bottomrule
\end{tabularx}
\end{table}

For $\Gamma=C_4=\langle g\rangle$, let $\chi(g)=\ii$.
The corresponding classes are listed in Table~\ref{tab:app-pointed-8-C4}.

\begin{table}[H]
\centering
\caption{Non-semisimple pointed classes of dimension $8$: $\Gamma=C_4$}
\label{tab:app-pointed-8-C4}
\begin{tabularx}{\textwidth}{@{}c X c@{}}
\toprule
Class & Datum $\DD_{\mathrm{pt}}$ or realization & Dual \\
\midrule
$\HH_8^{(12)}$ & $\kk C_4\otimes\Lambda_1=\AA(C_4,(1,1,0;1))$ & self-dual \\
$\HH_8^{(13)}$ & $\AA(C_4,(g,1,0;1))$ & $\HH_8^{(15)}$ \\
$\HH_8^{(14)}$ & $\AA(C_4,(g^2,1,0;1))$ & $\HH_8^{(16)}$ \\
$\HH_8^{(15)}$ & $\AA(C_4,(1,\chi,0;1))$ & $\HH_8^{(13)}$ \\
$\HH_8^{(16)}$ & $\AA(C_4,(1,\chi^2,0;1))$ & $\HH_8^{(14)}$ \\
$\HH_8^{(17)}$ & $\AA(C_4,(g^2,\chi^2,0;1))$ & self-dual \\
$\HH_8^{(18)}$ & $\AA(C_4,(g,1,1;1))$ & non-pointed \\
\bottomrule
\end{tabularx}
\end{table}

\subsection{Non-semisimple cases that are not pointed but have pointed duals}
\label{subsec:app-dual-pointed}
In dimensions at most $10$, there are exactly four non-trivial classes which are not pointed but whose duals are pointed, by~\cite[Theorems~5.1--5.4]{SSW24}.
They are listed in Table~\ref{tab:app-dual-pointed}.

\begin{table}[H]
\centering
\caption{Non-semisimple classes that are not pointed but have pointed duals}
\label{tab:app-dual-pointed}
\begin{tabularx}{\textwidth}{@{}c c c X@{}}
\toprule
Dimension & Class & Pointed dual & Additional property \\
\midrule
$6$ & $(\HH_6^{(4)})^*$ & $\HH_6^{(4)}$ & unique such class \\
$8$ & $(\HH_8^{(18)})^*$ & $\HH_8^{(18)}$ & does not have the Chevalley property \\
$8$ & $(\HH_8^{(19)})^*$ & $\HH_8^{(19)}$ & does not have the Chevalley property \\
$10$ & $(\HH_{10}^{(4)})^*$ & $\HH_{10}^{(4)}$ & unique such class \\
\bottomrule
\end{tabularx}
\end{table}

\subsection{Non-semisimple non-pointed cases with non-pointed duals}
For $\zeta\in\{\ii,-\ii\}$ and $\epsilon,\eta\in\{0,1\}$, let
$K_8(\zeta;\epsilon,\eta)$ be the superalgebra generated by $g,v,w,T$ with
\[
 |g|=|v|=\bar0,
 \qquad
 |w|=|T|=\bar1,
\]
and relations
\[
 g^2=1,
 \qquad gv=v,
 \qquad gw=-w,
 \qquad vw=wv=0,
\]
\[
 v^2=\frac12(1+g),
 \qquad
 w^2=\frac12(1-g),
\]
\[
 T^2=0,
 \qquad Tg=gT,
 \qquad Tv=(-1)^\epsilon vT,
 \qquad Tw=-(-1)^\epsilon wT.
\]
Its Hopf superalgebra structure is
\[
\begin{aligned}
\Delta(g)&=g\otimes g,
&\Delta(v)&=v\otimes v-\zeta w\otimes w,\\
\Delta(w)&=v\otimes w+w\otimes v,
&\Delta(T)&=g^\eta\otimes T+T\otimes1,\\
\varepsilon(g)&=\varepsilon(v)=1,
&\varepsilon(w)&=\varepsilon(T)=0,\\
S(g)&=g,
&S(v)&=v,
&S(w)&=\zeta w,
&S(T)&=-g^\eta T.
\end{aligned}
\]
The eight members of this family are pairwise non-isomorphic and
non-semisimple, and neither they nor their duals are pointed.
By~\cite[Section~5.5]{SW25}, they all have the Chevalley property and their bosonizations satisfy
\[
 \widehat{K_8(\zeta;\epsilon,\eta)}
 \cong H_{16}((-1)^{\epsilon+\eta}\zeta),
\]
where $H_{16}(\pm\ii)$ are the Hopf algebras of~\cite{CDMM04,GV10}.
For each fixed $\zeta$, the duality is
\[
 K_8(\zeta;0,0)^*\cong K_8(\zeta;0,0),
 \qquad
 K_8(\zeta;1,1)^*\cong K_8(\zeta;1,1),
\]
\[
 K_8(\zeta;0,1)^*\cong K_8(\zeta;1,0),
 \qquad
 K_8(\zeta;1,0)^*\cong K_8(\zeta;0,1).
\]
See \cite[Proposition~5.17 and Theorem~5.18]{SW25}.

\subsection{Semisimple cases with nonzero odd part}
Throughout this subsection, $p$ denotes an odd prime.

\subsubsection{Dimensions \texorpdfstring{$p^2$ and $2p$}{p2 and 2p}}
There is no non-trivial semisimple non-pointed Hopf superalgebra of dimension
$p^2$.
In dimension $2p$, there is exactly one, namely $(\kk D_{2p})^a$, and it is
self-dual.
For dimensions at most $10$, this gives the classes $(\kk D_6)^a$ and
$(\kk D_{10})^a$.
The first of these is the Hopf superalgebra $\AA_6$ of \cite{SW25}.

\subsubsection{Dimension \texorpdfstring{$4$}{4}}
Let $C_2=\langle x\mid x^2=1\rangle$ and let $\ppp(x)=1$.
The two classes are
\[
 \LL_4^{(1)}=\LL(C_2,(\id,\ppp;\lambdaa^+,1)),
 \qquad
 \LL_4^{(2)}=\LL(C_2,(\id,\ppp;\lambdaa^-,1)),
\]
where $\lambdaa^\pm(x)=\pm\ii$.
Both are self-dual.
In the notation of \cite{SW25}, these are $\AA_4(\ii)$ and
$\AA_4(-\ii)$, respectively.

\subsubsection{Dimension \texorpdfstring{$8$}{8}}
For $C_4=\langle x\mid x^4=1\rangle$, let $\ppp(x)=1$ and let $\iota$ denote inversion.
We use the tuple notation $(\lambdaa)_{C_4}$ introduced in Notation~\ref{not:cochain-values}.
The four classes with $\Gamma=C_4$ are listed in Table~\ref{tab:app-semisimple-C4}; in all four cases, $\muu=1$.

\begin{table}[H]
\centering
\caption{Semisimple eight-dimensional classes with $\Gamma=C_4$}
\label{tab:app-semisimple-C4}
\begin{tabular}{@{}ccccc@{}}
\toprule
$\HH$ & $\ttt$ & $(\lambdaa)_{C_4}$ & $\widehat{\HH}$ & $\HH^*$ \\
\midrule
$\LL_8^{(1)}$ & $\id$ & $\bigl(1,-1,-1,1\bigr)$ & $H_{c:\sigma_0}$ & $\LL_8^{(1)}$ \\
$\LL_8^{(2)}$ & $\id$ & $\bigl(1,-\ii,1,-\ii\bigr)$ & $H_{c:\sigma_1}$ & $\LL_8^{(10)}$ \\
$\LL_8^{(3)}$ & $\iota$ & $\bigl(1,\ii,1,\ii\bigr)$ & $H_{c:\sigma_1}$ & $\LL_8^{(12)}$ \\
$\LL_8^{(4)}$ & $\iota$ & $\bigl(1,-\ii,1,-\ii\bigr)$ & $H_{c:\sigma_1}$ & $\LL_8^{(11)}$ \\
\bottomrule
\end{tabular}
\end{table}

For $C_2^2=\langle y,z\mid y^2=z^2=1, yz=zy\rangle$, let $\ppp(y)=\ppp(z)=1$, and let $\vartheta$ be the automorphism interchanging $y$ and $z$.
We use the tuple notation $(\lambdaa)_{C_2^2}$ introduced in Notation~\ref{not:cochain-values}.
The normalized non-trivial $2$-cocycle used throughout is
\[
 \muu_0(y^rz^s,y^uz^v)=(-1)^{r(u+v)}.
\]
Table~\ref{tab:app-mu0} gives all its values; the row index is the first argument and the column index is the second.

\begin{table}[H]
\centering
\caption{Values of the normalized cocycle $\muu_0$ on $C_2^2$}
\label{tab:app-mu0}
\begin{tabular}{@{}c rrrr@{}}
\toprule
$\muu_0(g,h)$ & $h=1$ & $h=y$ & $h=z$ & $h=yz$\\
\midrule
$g=1$  & $1$& $1$& $1$& $1$\\
$g=y$  & $1$&$-1$&$-1$& $1$\\
$g=z$  & $1$& $1$& $1$& $1$\\
$g=yz$ & $1$&$-1$&$-1$& $1$\\
\bottomrule
\end{tabular}
\end{table}

The remaining eight classes are listed in Table~\ref{tab:app-semisimple-C22}.

\begin{table}[H]
\centering
\caption{Semisimple eight-dimensional classes with $\Gamma=C_2^2$}
\label{tab:app-semisimple-C22}
\begin{tabular}{@{}cccccc@{}}
\toprule
$\HH$ & $\ttt$ & $(\lambdaa)_{C_2^2}$ & $\muu$ & $\widehat{\HH}$ & $\HH^*$ \\
\midrule
$\LL_8^{(5)}$ & $\vartheta$ & $\bigl(1,1,1,1\bigr)$ & $\muu_0$ & $H_{d:1,1}$ & $\LL_8^{(5)}$ \\
$\LL_8^{(6)}$ & $\vartheta$ & $\bigl(1,\ii,\ii,1\bigr)$ & $1$ & $H_{d:-1,1}$ & $\LL_8^{(9)}$ \\
$\LL_8^{(7)}$ & $\id$ & $\bigl(1,\ii,\ii,1\bigr)$ & $1$ & $H_{d:-1,1}$ & $\LL_8^{(7)}$ \\
$\LL_8^{(8)}$ & $\id$ & $\bigl(1,-\ii,-\ii,1\bigr)$ & $1$ & $H_{d:-1,1}$ & $\LL_8^{(8)}$ \\
$\LL_8^{(9)}$ & $\id$ & $\bigl(1,\ii,\ii,1\bigr)$ & $\muu_0$ & $H_{d:-1,1}$ & $\LL_8^{(6)}$ \\
$\LL_8^{(10)}$ & $\id$ & $\bigl(1,-\ii,\ii,-1\bigr)$ & $1$ & $H_{d:-1,-1}$ & $\LL_8^{(2)}$ \\
$\LL_8^{(11)}$ & $\id$ & $\bigl(1,\ii,-\ii,-1\bigr)$ & $\muu_0$ & $H_{d:-1,-1}$ & $\LL_8^{(4)}$ \\
$\LL_8^{(12)}$ & $\id$ & $\bigl(1,-\ii,\ii,-1\bigr)$ & $\muu_0$ & $H_{d:-1,-1}$ & $\LL_8^{(3)}$ \\
\bottomrule
\end{tabular}
\end{table}

By Propositions~\ref{prp:K5-D8},~\ref{prp:K6-KP}, and~\ref{prp:Masuoka-split}, we have
\[
\LL_8^{(5)}\cong(\kk D_8)^a,
\qquad
\LL_8^{(6)}\cong(\AKP)^s,
\qquad
\LL_8^{(7)}\cong\LL_4^{(1)}\otimes\kk C_2,
\qquad
\LL_8^{(8)}\cong\LL_4^{(2)}\otimes\kk C_2.
\]

\subsection*{Acknowledgment}
The first author (T.S.) was supported by JSPS KAKENHI Grant Numbers JP22K13905 and JP26K06774.
The second author (K.S.) was supported by JSPS KAKENHI Grant Number JP24K06676.
The authors acknowledge the assistance of ChatGPT (GPT-6) in proofreading the manuscript and improving the clarity of the exposition.
The authors take full responsibility for the mathematical content and the final text.

\providecommand{\href}[2]{#2}

\end{document}